\input graphicx

\vsize=9truein

\def \bbu {{\hskip -.03in}}

\def \ST {{\cal ST}}
\def \tM {{ \widetilde M}}

\def \BZ  {{\bf Z}}

\def\multi#1{\vbox{\baselineskip=0pt\halign{\hfil$\scriptstyle\vphantom{(_)}##$\hfil\cr#1\crcr}}}

\def \BZ {{\bf Z}}

\def\tttt #1{{\textstyle{#1} }}
\def\ul #1{{\underline{#1} }}
\def\ol #1{{\overline{#1} }}

\def \LAr {\leftarrow}

\def \CT {{\cal T}}

\def \DA {\downarrow}

\def \magstep#1 {\ifcase#1 1000\or 1200\or 1440\or 1728\or 2074\or 2488\fi\relax}

\def\la{{\lambda}}

\def \CA {{\cal A}}

\overfullrule=0pt
\baselineskip 14pt
\settabs 8 \columns
\+ \hfill&\hfill \hfill & \hfill \hfill   & \hfill \hfill  & \hfill \hfill   &
 & \hfill \hfill\cr
\parindent=.5truein
\hfuzz=2.44182pt
\hsize 6.5truein
\font\ita=cmssi11  
\font\small=cmr6
\font\title=cmbx11 scaled\magstep2
\font\normal=cmr11 
\font\small=cmr6

\font\bol=cmbx11

\def \con {\subset}

\def \ggg {\gamma}

\def\PI{\Pi}
\def \-> {\rightarrow}
\def\LL{\big\langle}

\def\RR {\big\rangle}

\def \ue {{\underline e}_1}

\def\DD {\Delta}

\def\OM {\Omega}

\def\om {\omega}
\def\la {\lambda}
\def\La {\Lambda}
\def \RA {\rightarrow}
\def \LA {\lefttarrow}
\def\xon {x_1,x_2,\ldots ,x_n}

\def \sas {\vskip .06truein}
\def\sa{{\vskip .125truein}}

\def\sap{{\vskip .25truein}}

\def \eee {\epsilon}

\def\aaa {\alpha}
\def\bbb {\beta}
\def\ddd{\delta}

\def\ggg {\gamma}

\def\noc {\supset}
\def \ses {\enskip = \enskip}
\def \sps {\, + \,}

\def \sscs {\, ,}

\def \sms {\, - \,}

\def \scs {\, , \,}
\def \ess {\enskip}

\def \ssp {\hskip .25em}
\def \bigsp {\hskip .5truein}
\def \part {\vdash}

\def \DD {\Delta}

\normal

\vsize=9truein
\sap
\def\today{\ifcase\month\or
January\or February\or March\or April\or may\or June\or
July\or August\or September\or October\or November\or
December\fi
\space\number\day, \number\year}

\def \OM {{\Omega}}

\def \RA {{ \rightarrow }}
\def \LA {{ \leftarrow }}

\def \CS {{\cal S}}

\def \BQ {{\bf Q}}

\def \II {1}
\def \om {\omega}

\def \TH {{\widetilde H}}

\def\xon {x_1,x_2,\ldots ,x_n}

\def \BQ {{\bf Q}}
\def \om {\omega}

\font\small=cmr6
\def \scs {\ssp , \ssp}
\def \ess {\enskip}
\def \ssp {\hskip .25em}
\def \bigsp {\hskip .5truein}
\def \part {\vdash}
\font\title=cmbx11 scaled\magstep2
\font\normal=cmr11 
\def\today{\ifcase\month\or
January\or February\or March\or April\or may\or June\or
July\or August\or September\or October\or November\or
December\fi
\space\number\day, \number\year}
\null

\def \II {{\rm I}}

\def \I {{\rm I}}

\def \TM {\widetilde M}

\def \ck {c_{\mu,\nu}^{(k)}} 
\def \dk {d_{\mu,\nu}^{(k)}} 

\def \and {\ess\ess\ess\hbox{and}\ess\ess}

\headline={  \small \today \hfill {\small  Stanley-Macdonald Pieri Rules}\hfill$\ess\ess$  \folio}
 \footline={\hfill }
 
 \null
 \vskip -.75in
 
 \centerline{\title  Some new applications}
  \centerline{\title of the Stanley-Macdonald Pieri Rules}
\sas

  \centerline{\bol  A.M. Garsia, J. Haglund, G. Xin and M. Zabrocki}
\sas

\centerline{\ita Dedicated to Richard   Stanley  
for his 70th birthday}

\sas

\noindent{\bol Abstract}

In a seminal  paper [22] Richard Stanley derived Pieri rules for the Jack symmetric function basis. These rules were   extended   by Macdonald to his now famous symmetric function basis. The original form  of these rules had a forbidding complexity that made them difficult to use in explicit calculations.
In the early 90's it was discovered [10] that, due to massive cancellations, the 
dual rule, which  expresses  skewing by $e_1$ the modified Macdonald polynomial 
$\TH_\mu [X;q,t]$, can be given a very simple combinatorial form in terms of corner weights of the Ferrers' diagram of $\mu$.  A similar formula was later obtained by the last named author for the multiplication of $\TH_\mu [X;q,t]$ 
by $e_1$, but never published. In the years that followed we have seen some truly remarkable uses of these two Pieri rules in establishing highly non trivial combinatorial results 
in the Theory of Macdonald polynomials. This theory has recently been
spectacularly enriched by various Algebraic Geometrical results in   the  works of  Hikita [15], 
Schiffmann [19], Schiffmann-Vasserot  [20],  [21], A. Negut [18] and Gorsky-Negut [12].
This development opens up the challenging task of deriving their results 
by purely Algebraic Combinatorial methods. Substantial progress in this task was obtained in [3]. In this paper we present the progress obtained by means of
 Pieri rules.
  \sa
  
  \noindent{\bol I. Introduction}
  
  In this paper, as in [3], the main actors are the family of symmetric function operators $D_k$ and $D_k^*$ introduced in [6], whose action on a symmetric function $F[X]$ may be written on the form
$$
 \eqalign{
a) \ess\ess\ess D_k F[X]
 &\ses F[X+\tttt{M\over z}]\sum_{i\ge 0}(-z)^ie_i[X]\Big|_{z^k}
 \bigsp(\hbox{with $M=(1-t)(1-q)$})
 \cr
b) \ess\ess\ess    D_k^*F[X]
&\ses F[X-\tttt{\TM\over z}] \sum_{i\ge 0}z^ih_i[X]\Big|_{z^k}
\bigsp \big(\hbox{$\TM=(1-1/t)(1-1/q)$}\big)
\cr}
\eqno \I.1
$$
where expressions are given in plethystic notation which we shall review in Section 1.

We will focus here, as in [3], on the algebra $\CA$ of symmetric function operators generated by the family  $\{ D_k\}_{k\ge 0}$.  It was shown in [3] that $\CA$  is bi-graded by assigning the generator $D_k$ bi--degree $(1,k)$ .  Its connection to the Algebraic Geometrical developments is that $\CA$ is a concrete realization of a portion of the  Elliptic Hall Algebra  studied by Schiffmann and Vasserot in [19] ,[20] and [21].  

 In particular $\CA$  contains a distinguished family of operators  $\{Q_{u,v}\}_{u,v\ge0}$ of bi-degree given by their index that can be shown to play a central role in the connection between Macdonald Polynomials  and Parking Functions  (see [4]). For nonnegative, and co-prime index pair $(m,n)$ the   construction of the  operators 
 $Q_{m,n}$ is quite simple. We write 
$$
Split(m,n)= (a,b)+(c,d)
\ess\ess\ess 
\hbox{if and only if }\ess\ess\ess  
\cases {a) & $(m,n)= (a,b)+(c,d)$\cr
b)  & $det\Big\| \matrix{a & c\cr b & d\cr}\Big\| =1$\cr}
.
\eqno \I.2
$$ 
Geometrically this simply says  that $(a,b)$ is the lattice point closest to
the segment $(0,0)\RA (m,n)$ and $(0,0) \RA (a,b) \RA(m,n)$ is the counter-clockwise order of the vertices of a non-trivial  triangle. These conditions force also  the pairs
$(a,b)$ and $(c,d)$ to be co-prime and we  can  recursively define
$$
Q_{m,n}\ses \tttt{1\over M} [Q_{c,d},Q_{a,b}] = 
\tttt{1\over M} (Q_{c,d} Q_{a,b} - Q_{a,b} Q_{c,d})
\eqno \I.3
$$
with base cases
$$
a)\ess\ess Q_{1,0}=D_0\ess\ess\ess
\hbox{and}\ess\ess\ess
b)\ess\ess Q_{0,1}=- \ul e_1
\eqno \I.4
$$
where  for a symmetric function $f$ we let ``$\ul f$''
denote the  operator, ``{\ita multiplication by $f$}.''

Now from the  definition in I.1 a) it easily follows that 
$$
D_k\ses \tttt{1\over M}[D_{k-1},\ul e_1]\bigsp( \hbox{for all $k\ge 1$})
\eqno \I.5
$$
and I.3 then yields
$$
Q_{1,k}= D_k.
\eqno \I.6
$$
Thus for $m,n$ positive integers we may replace the recursion in I.3 by
$$
Q_{m,n}\ses \cases 
{ \tttt{1\over M} [Q_{c,d},Q_{a,b}] & if $m>1$, \cr\cr
D_n & if $m=1$.
}
\eqno \I.7
$$
This implies that  $Q_{m,n}\in \CA$ for all  $(m,n)$ positive and co-prime. 
 
The definition of the operators   $Q_{u,v}$ for non co-prime pairs $(u,v)$ and the proof of some of their remarkable properties
was carried out in [3] by means of two basic tools.  

To be more explicit we need some auxiliary material.
Firstly, we will  write $(u,v)=(km.kn)$ with $(m,n)$ a co-prime pair and $k=gcd(u,v)$. Next note that 
from $Split(m,n)=(a,b)+(c,d)$ we derive that
$$
(u,v)\ses (im+a,in+b)+(jm+c,jn+d)\ess\ess\ess\ess
(\hbox{for all $i+j=k-1$} )
\eqno \I.8
$$
This given,  in [3] it is shown that we may set 
$$
Q_{u,v}\ses \tttt{1\over M} [Q_{jm+c,jn+d},Q_{im+a,in+b}]
\ess\ess\ess\ess
(\hbox{for any $i+j=k-1$} )
\eqno \I.9
$$
by proving that all the  operators on the right hand side  are the same. 
Since both  pairs $(im+a,in+b)$ and $(jm+c,jn+d)$ always turn out to be  co-prime, this  shows that 
$Q_{u,v}$ is well defined,
  and also  shows that  
$Q_{u,v}\in \CA $ for any  $u,v\ge 1$.

The equality of the operators occurring on the right hand side of I.9 was obtained in [3] by two steps. 
In the first step I.9 is established for
$(m,n)=(1,1)$ and in the second step I.9 is derived in full generality 
from this special case by means 
of a natural action of $SL_2[\BZ]$
on $\CA$ 
which preserves  all the identities  satisfied by the operators $D_k$.
In particular a proof is given in  [3] that for the generators  
$\Big[\matrix{a & c\cr b & d}\Big]=\Big[\matrix{1 & 1\cr 0 & 1}\Big]$ 
and $ \Big[\matrix{a & c\cr b & d}\Big]=\Big[\matrix{1 & 0\cr1 & 1}\Big]$
we have
$$
\Big[\matrix{a & c\cr b & d}\Big] Q_{m,n}= Q_{am+cn,bm+dn}.
$$
This given, it is first shown in [3] that
the case $(m,n)=(1,1)$ of I.9 is 
a consequence of the
commutator  identity
$$
(D_aD_b^*-D_b^*D_a)\, P[X]\ses 
M\tttt{(qt)^{-a}\over tq-1} 
h_{a+b}\big[X(1-tq)\big] P[X]  
\ess\ess\ess\ess
\hbox{(for $a+b>0$)}
\eqno  \II.10
$$
In fact, by setting
$$
D_k=Q_{1,k}  \ess \ess\ess \ess \hbox{and}\ess \ess\ess \ess
D_k^*= -(qt)^{k-1} Q_{-1,k}
\eqno \II.11
$$
I.10 becomes the operator identity
$$
\tttt{1\over M}[Q_{-1,b}\scs Q_{1,a}]
\ses \tttt{qt\over qt-1}
\, \ul h_{a+b}[X(1/qt-1)].
\eqno \II.12
$$
 Next  it is shown in [3] that for all co-prime $(m,n)$ we have
$$
\nabla Q_{m,n}\nabla^{-1}\ses  \Big[\matrix{1 & 1\cr 0 & 1}\Big] 
 Q_{m,n}
\ses
Q_{m+n,n}
\eqno \II.13
$$
where $\nabla$ is the operator,
introduced  in [1] with eigen-functions the  modified Macdonald basis 
$\{\TH_\mu[X;q,t]\}_\mu$.
Thus conjugating I.12 by $\nabla$, gives
$$
\tttt{1\over M}[Q_{b-1,b}\scs Q_{a+1,a}]
\ses \tttt{qt\over qt-1}
\, \nabla \ul h_{a+b}[X(1/qt-1)]\nabla^{-1}
\ess\ess\ess\ess
\hbox{(for $a+b>0$)}.
\eqno \II.14
$$
Since for $(m,n)=(1,1)$ we have
$(a,b)=(1,0)$ and $(c,d)=(0,1)$
we see that I.14  implies
that for $(m,n)=(1,1)$ all the operators on the right hand side of I.9 are identical  and 
we may thus set 
$$
Q_{k,k}\ses  \tttt{qt\over qt-1}
\, \nabla \ul h_{k}[X(1/qt-1)]\nabla^{-1}
.
\eqno \II.15
$$
The final step is obtained by showing that the  identities in I.9 are 
simply images of I.14 and I.15  by the   $SL_2[Z]$ action.
\sap

In this paper, using the symmetric function tools created in the 90's 
in the study of Macdonald polynomials,  most particularly  in [8], [6], [2], 
we develop  a parallel variety of identities  by letting the operators $Q_{k,-1} $ and $Q_{k,1}$ play the role that $Q_{1,k}$ and $Q_{-1,k}$ play in [3]. 
The  fact that for 
$$
U=
\Big[\matrix{0 & 1\cr -1 & 0}\Big]\ses
\Big[\matrix{1 & 1\cr 0 & 1}\Big]
\Big[\matrix{1 & 0\cr 1 & 1}\Big]^{-1}
\Big[\matrix{1 & 1\cr 0 & 1}\Big]
$$ 
we have
$$
UQ_{1,k}=Q_{k,-1} \scs\ess\ess\ess\ess
UQ_{-1,k}=Q_{k,1}
$$
guided us to a number of surprising discoveries. In particular it turns out that an identity that was discovered in 2008
in the research that yielded the results in [11] may be viewed as an image by 
$U$ of the identity in I.12. 

To state our results we need to recall some notational conventions. To begin, we will identify partitions with their 
French Ferrers diagrams. Next, for a cell $c\in \mu$ we let $l_\mu(c)$,
$a_\mu(c)$, $l'_\mu(c)$ and $a'_\mu(c)$ denote the ``{\ita leg}'', ``{\ita arm}'', ``{\ita coleg}' and \hbox{``{\ita coarm }''}, of $c$ in $\mu$ (as in [17]). 
Because we are using French notation,  these parameters count the number of cells of $\mu$ that are  respectively strictly North, East, South and West of $c$ in $\mu$. We  then set
$$
B_\mu(q,t)= \sum_{c\in \mu}t^{l'_\mu(c)}q^{a'_\mu(c)},\bigsp 
\Pi_\mu(q,t)= \prod_{c\in \mu}\bbu \ ^{(0,0)}(1-t^{l'_\mu(c)}q^{a'_\mu(c)})
\eqno \II.16
$$
where the superscript $(0,0)$ in the product is to  avoid the vanishing factor. In addition we set
$$
T_\mu =\prod_{c\in \mu }t^{l'_\mu(c)}t^{a'_\mu(c)}
,\bigsp w_\mu(q,t)=\prod_{c\in \mu}(q^{a _\mu(c)} -t^{l _\mu(c)+1})(t^{l _\mu(c)} -q^{a _\mu(c)+1}),
\eqno \II.17
$$
The modified Macdonald polynomials  $\{\TH_\mu[X;q,t]\}_\mu$ we work with here form the unique symmetric function basis that is (in dominance order)  upper triangularly   related to the modified Schur basis $\{s_\la\big[\tttt{X\over t-1}\big]\}_\la$ and satisfies the 
orthogonality condition
\vskip -.15in
 $$ 
\bigsp \LL \TH_\la\scs \TH_\mu\RR_*\ses \chi(\la=\mu) w_\mu(q,t),
 \eqno \II.18
 $$
 \noindent
where $\LL\scs \RR_*$ is the deformation   
\vskip -.1in
$$
\LL p_\la\scs p_\mu \RR_*\ses 
(-1)^{|\mu|-l(\mu)}   \prod_{i=1}^{l(\mu)} (1-t^{\mu_i})(1-q^{\mu_i})\ssp z_\mu\ssp
\chi(\la=\mu) 
\eqno \II.19
$$
\vskip -.05in
 \noindent
 of the Hall  scalar product
$\LL p_\la\scs p_\mu \RR=z_\mu\ssp
\chi(\la=\mu)$.
 We will call $\LL\scs \RR_*$  the ``{\ita star scalar product}''.

Let us also recall   that the operator $\nabla$ is defined by setting 
$$
\nabla \TH_\mu[X;q,t]\ses T_\mu \TH_\mu[X;q,t].
\eqno \II.20
$$
It is also shown in [6] that
$$
D_0\TH_\mu[X;q,t]=\big(1-MB_\mu(q,t)\big)\TH_\mu[X;q,t].
\eqno \II.21
$$
Both $\nabla$ and $D_0$ are special cases  of a commuting family of operators defined in [6] by setting for a symmetric function $F$ 
$$
\DD_F \TH_\mu[X;q,t]\ses F[B_\mu(q,t)] \TH_\mu[X;q,t].
\eqno \II.22
$$ 
Finally, the family of operators $ Q_{m,n}$ with $m,n$ co-prime is extended
to the fourth lattice quadrant by setting
$$
Q_{m,-n}\ses Q_{m,n}^{\perp^*}
\eqno \II.23
$$
where the symbol ``$\perp^*$'' denotes the operation of taking the 
adjoint of an operator with respect to the star scalar product.
\sap

This given, our first result can be stated as follows.
\sas

\noindent{\bol Theorem I.1}
 
{\ita For all $m\ge 1$ we have}

$$
Q_{m,0}\ses  \tttt{qt \over qt-1}\DD_{h_m \big[(MX-1)(1/qt-1) \big]}
\eqno \II.24
$$

Notice that, by I.21, this identity  may be viewed as the extension to $m>1$ of the equality $Q_{1,0}=D_0$.
Notice further that since $Q_{0,m}= {qt\over  qt-1}\ul
h_m[X(1/qt-1)]$ and the collection 
$\big\{\prod_{i=1}^{\l(\la)}Q_{\la_i,0}\big\}_\la$ may be taken as a basis for the family of symmetric function multiplication operators, we can derive 
from Theorem I.1 and the identity $U Q_{0,m}=Q_{m,0}$
   the following  truly remarkable result.
\sas

\noindent{\bol Theorem I.2}

{\ita The action of the $2\times 2$ matrix $U$ on a symmetric function operator $\ul F $ may be expressed by the identity}
$$
U\, \ul F[X]\ses \Delta_{F[(MX-1)]}
\eqno \II.25
$$

Another significant fact that emerges from our findings here is that 
while viewing the operators $Q_{m,n}$ as non-commutative polynomials in the family  $Q_{1,n}=\Big[\matrix{0 & 1\cr 1 & 1}\Big]^nQ_{0,1}=D_n$  expresses their action by constant term formulas, our present way of viewing the operators $Q_{m,n}$ as non-commutative polynomials in the family  $Q_{m,1}=\Big[\matrix{1 & 1\cr 0 & 1}\Big]^mQ_{0,1}=\nabla^mQ_{0,1}\nabla^{-m}$ expresses their action by standard tableaux expansions. 

To be more precise  we need notation.
Let ${\cal ST}_n$ be the set of all standard tableaux with labels $1,2,\dots,n$ and $ST(\mu)$ be all the standard tableaux
of shape $\mu$. 
For a given $T\in {\ST}_n$, we set $w_T(k)=q^{j-1}t^{i-1}$, if the label $k$ of $T$ is in the $i$-th row and the $j$-th column. 
This given, our simplest result in this context  may be stated as follows.

\sas

\noindent{\bol Theorem I.3}

 {\ita The operators $Q_{u,v}$ may be written as a linear combination of the family    
 $
\big\{Q_{a_n,1}\cdots Q_{a_2,1}Q_{a_1,1}\big\}_{a_i\ge 0},
 $
whose action on $(-1)^n$ may be expressed in the form
$$
 {1\over M^n}Q_{a_n,1}\cdots Q_{a_2,1}Q_{a_1,1} (-1)^n\ses
\sum_{\mu\part n}{\TH_\mu[X;q,t]\over w_\mu}\bbu
\sum_{T\in ST(\mu)}\bbu \prod_{k=2}^n w_T(k)^{a_k }
 \tttt{ 1-w_T(k) qt   \over 1-qt }
\bbu \prod_{1\le h<k \le n}\bbu
f\big(\tttt{  w_T(h)\over w_T(k) }\big)
\eqno \II.26
$$
where for convenience we have set $f(u)= {(1-u/t)(1-u/q)\over (1-u)(1-u/qt)}$.
}
\sas

We should mention that a similar tableau expansion was first obtained in [9] and  in fact our proof of I.26 follows very closely the arguments used in [9].
This is precisely where the Stanley-Macdonald Pieri rules play a crucial role. We must also mention that an entirely analogous standard tableaux expansion 
for this particular  action of the $Q_{m,n}$ operators (for $(m,n)$ co-prime) was given in [12].
Nevertheless, a distinguishing feature of I.26 is that it is a closed formula.
By contrast, the expansions in  [12] as well as in [9] are really algorithms,
in that the substitutions $x_k\RA1/w_T(k)$ used there are 
to be carried out iteratively (as we will see in section 3) due to the fact
that the kernel to which these substitutions are performed has denominators that vanish under these substitutions.

Here the crucial tool that makes all this possible is the following identity which may be viewed as the image by $U$ of the identity in I.12.
\sas

\noindent{\bol Theorem I.4}

{\ita 
For all $a,b\ge 0$ we have}
$$
\tttt{1\over M}[Q_{a,1}\scs Q_{b,-1}]
\ses \tttt{qt\over qt-1}
\DD_{ h_{a+b}[(MX-1)(1/qt-1)]}.
\eqno \II.27
$$

\sa

Our presentation is divided into  3 sections. In the first section we will review some of the identities and  definitions that we need in the present development which were introduced or proved elsewhere. In the second section we give the proof of Theorem  I.1  and complete  
our treatment of the algebra $\CA$ as generated by the operators
$Q_{m,1}$. In the third section we give the heretofore unpublished combinatorial  argument that
derives, from the Stanley-Macdonald Pieri rules,    
``{\ita corner weights}'' expressions for the coefficients $d_{\mu,\nu}$ in the expansion
$$
e_1\TH_\nu[X;q,t]\ses\sum_{\mu\LA \nu}d_{\mu\nu}\TH_\nu[X;q,t].
\eqno \II.28
$$
This done, in this  section we prove  two versions of our standard tableaux expansion formulas one of which is  I.26. We terminate this section and the paper by pointing out that a recent result of Bergeron-Haiman [5]  shows that the Pieri coefficients  $d_{\mu\nu}$ and  $c_{\mu,\nu}$ are not as limited as may appear on the surface.
  \sa

\noindent {\bol 1. Preliminaries}

\sas

The space of symmetric polynomials  with coefficients in $\BQ[q,t]$  will be denoted $\Lambda$. The subspace of
homogeneous symmetric polynomials of degree $m$ will be denoted
$\Lambda^{=m}$.
We will seldom work with symmetric polynomials expressed in terms of variables but rather express them
in terms of one of the classical symmetric function bases
$\{m_\la\}_\la$, $\{p_\la\}_\la$, $\{h_\la\}_\la$,
$\{e_\la\}_\la$ and $\{s_\la\}_\la$ ({\it Schur}).

We recall that the fundamental involution $\om$ may be defined by setting 
for the power basis
$$
\om p_\la\ses (-1)^{n-k}p_\la\ses  (-1)^{|\la |-l(\la) }p_\la
\eqno 1.1$$
where for any vector $\ess v=(v_1,v_2,\ldots,v_k)$ we set
 $
\ess |v|=\sum_{i=1}^k v_i 
$
and 
$ l(v)=k $.

\sas

In dealing with symmetric function identities, specially with those arising
in the Theory of Macdonald Polynomials, we find it convenient
and often indispensable to use 
plethystic notation. This device has a straightforward definition which can be
verbatim implemented in MAPLE or MATHEMATICA
for computer experimentation. We simply set for
any expression $E=E(t_1,t_2 ,\ldots )$ and any power symmetric function $p_k$
$$
p_k[E]\ses E( t_1^k,t_2^k,\ldots ).
\eqno 1.2
$$
This given, for any symmetric function $F$ we set
$$
F[E]\ses Q_F(p_1,p_2, \ldots )\Big|_{p_k\RA E(t_1^k,t_2^k,\ldots )}
\eqno 1.3
$$
where $Q_F$ is the polynomial yielding the expansion of $F$ in 
terms of the power basis.

A paradoxical but necessary property of plethystic substitutions is 
that 1.3 requires  
$$
p_k[-E]\ses -p_k[E].
\eqno 1.4
$$
This notwithstanding, we will still need to carry out ordinary  changes 
of signs. To distinguish it from the ``{\ita plethystic}''  minus sign,
we will carry out the \hbox{``{\ita ordinary}''} sign change   
  by means
of a new variable  $\eee$ 
which outside of the plethystic bracket  is simply replaced by $-1$. 
For instance, these conventions give for $X_k=x_1+x_2+\cdots +x_n$
$$
p_k[ -\eee X_n]\ses  -\eee^k\sum_{i=1}^nx_i^k\ses  (-1)^{k-1} \sum_{i=1}^nx_i^k
$$  
Thus  for any symmetric function $F\in \Lambda$ and any expression $E$ we have
 $$ 
\om F[E]\ses F[-\eee E]
\eqno 1.5
$$ 
In particular,  if $F\in \Lambda^{=k} $  we may also rewrite this as
$$
F[-E]\ses    (-1)^k \om F[ E].
\eqno 1.6
$$
The formal power series
$$
\OM\ses exp\Big(\sum_{k\ge 1}{p_k\over k}\Big)
$$
 combined with plethysic substitutions    provides   
a powerful way of dealing with the many generating functions occurring in
our manipulations. In fact,  for  any given expression $E$ we will set
$$
\OM[E]\ses exp\Big(\sum_{k\ge 1}{p_k[E]\over k}\Big)
$$
and since for any two expressions $A,B$ 1.3 gives
$$
p_k[A+B]\ses p_k[A]+p_k[B]
\eqno 1.7
$$
We derive from this the fundamental  formula
$$
\OM[A+B]\ses\OM[A ]\, \OM[ B] 
\eqno 1.8
$$
In particular for  $A=\sum_{i=1}^n a_i $ and $B=\sum_{j=1}^m b_j $ we also get
$$
\OM[ z(A- B)]= {\prod_{j=1}^m  (1- b_jz) 
\over \prod_{i=1}^n( 1- a_iz)   }
\eqno 1.9
$$Clearly,  for any two expressions $A,B$  
we can   view  $\OM[z(A-B)]$   as the generating functions of the homogeneous 
symmetric functions plethystically evaluated at $A-B$
 $$
\OM[z(A-B)]= \sum_{m\ge 1} z^m h_m[A-B] 
$$
In particular, by equating coefficients of $z^m$ on both sides of 1.9,    we get 
(using 1.6)
$$
h_m[A-B]\ses \sum_{r=0}^m h_{m-r}[A]h_r[-B]\ses  \sum_{r=0}^m h_{m-r}[A](-1)^re_r[B]
$$
In particular it follows from this that
$$
h_m[(1-t)(1-q)]\ses
\cases{
 (1-t)(1-q){ \sum_{i=0}^{m-1}} (qt)^i &if $m>0$\cr\cr
 1 & if $m=0$\cr
 } 
\eqno 1.10
$$

\sas

The following  facts  (proved in [6]) will play a basic role here 

 \noindent{\bol Proposition 1.1} 

{\ita $D_k$ and $D_k^*$ are $*$-adjoint to $(-1)^k D_{-k}$ and
$(-qt)^kD_{-k}^*$ respectively. Moreover they are related to 
 the modified Macdonald polynomials $\TH_\mu[X;q,t]$
 and $\nabla$ by the  identities
$$
\matrix{
&(i )& \ess\ess D_0\, \TH_\mu= -D_\mu(q,t) \, \TH_\mu
\ess\ess 
&(i )^*&
D_0^*\, \TH_\mu= -D_\mu(1/q,1/t) \, \TH_\mu
\cr
&(ii)&\ess\ess  D_k\, \ue - \ue\, D_k = M\, D_{k+1} 
&(ii )^*&
D_k^*\, \ue - \ue\, D_k^* = \, - \widetilde M\, D_{k+1}^* 
\cr
&(iii)&\ess\ess\nabla \, \ue \nabla^{-1}= -D_1  
&(iii )^*&
\hskip -.06in\nabla \, D_1^* \nabla^{-1}= \ue
\cr
&(iv)&\ess\ess \nabla^{-1} \, e_1^\perp \nabla =  {\textstyle {1\over M}}D_{-1}   
&(iv )^*&
\ess\ess \nabla^{-1} \, D_{-1}^* \nabla = -{\widetilde M\,} e_1^\perp
\cr
&(v)&\ess\ess
D_ke_1^\perp\sms e_1^\perp D_k=D_{k-1}
&(v )^*&
\ess\ess 
D_k^*e_1^\perp\sms e_1^\perp D_k^*= -D_{k-1}^*
\cr
}
\eqno 1.11
$$
 with $e_1^\perp$ the Hall  scalar product adjoint of multiplication by $e_1$, $\widetilde M=(1-{1/ t})(1-{1/q})$  and 
}
$$
D_\mu(q,t)= MB_\mu(q,t)-1
\eqno 1.12
$$
As in [3] our starting point are  the identifications
$$
a)\ess\ess Q_{0,1}=-\ul e_1\scs
\bigsp\bigsp
b)\ess\ess Q_{1,0}= D_0\scs
\eqno 1.13
$$
Thus it follows from 1.11  $(ii)$, $(iii)$ and the  definition in I.3 that
$$
a)\ess\ess Q_{1,k}=D_k
\bigsp\bigsp
b)\ess\ess Q_{1,1}=\nabla Q_{0,1}\nabla^{-1}
\eqno 1.14
$$
Now it is shown in [3] that the definition in I.3 combined with 1.13 b) 
implies the following fundamental identity.
\sas

\noindent{\bol Proposition 1.2}

{\ita For any co-prime pair $m,n$ we have}
$$
Q_{m+n,n}\ses \nabla Q_{m,n} \nabla^{-1}.
$$
In particular it follows that we also have
$$
Q_{m,1}= \nabla^m Q_{0,1}\nabla^{-m}= - \nabla^m e_1\nabla^{-m} .
\eqno 1.15
$$
For notational  convenience, here and after we may use the symbol
``$T_m $'' to represent  ``$- \nabla^m e_1\nabla^{-m}$''.
The following basic fact is also an immediate consequence of the definition in I.3
\sas

\noindent{\bol Theorem 1.1}

{\ita For any co-prime   pair $m,n$ we may set }
$$
Q_{m,n}\ses 
\cases {
\tttt{1\over M}[Q_{c,d},Q_{a,b}] & if $n>1$
 and $Split(m,n)=(a,b)+(c,d)$\cr \cr
T_m & if $n=1$.
\cr
}
\eqno 1.16
$$
\vskip -.1in

\noindent
For instance, since we have (see adjacent figure)
 
 \vskip  -.01in
 \noindent
 $
 \ess\ess\ess \includegraphics[scale=0.26]{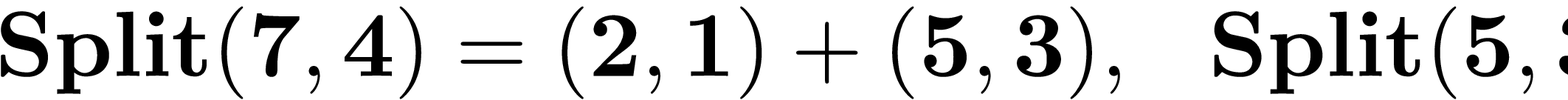}  
  $
  
\vskip-.6in
 \hfill
 $
  \includegraphics[width=1.2 in]{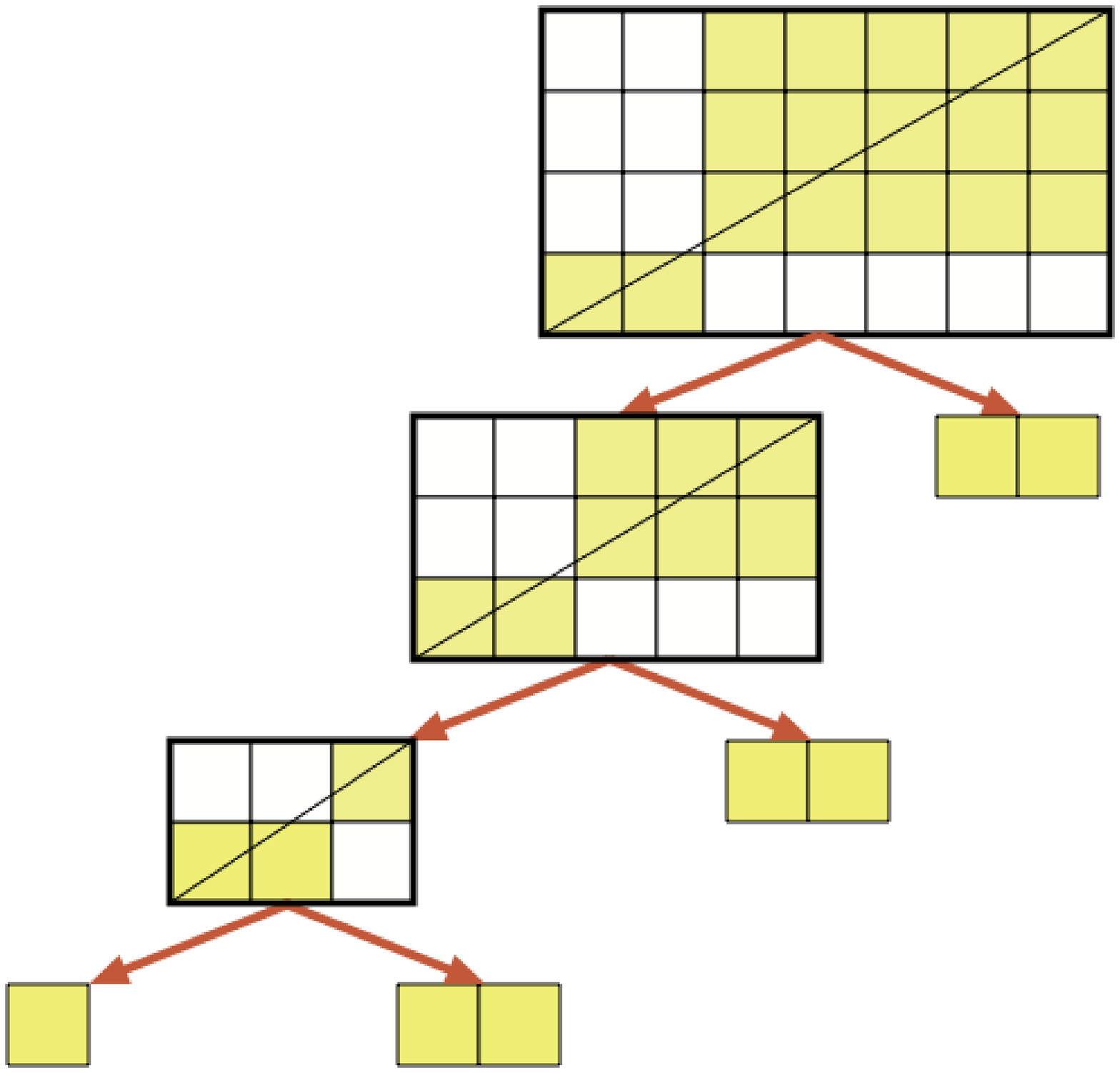} 
 $
 
 \vskip-.6in
 \noindent
The definition in 1.16 gives
\vskip-.1in
\noindent

$$
Q_{7,4}=\tttt{1\over M}[Q_{5,3},T_2]\scs \ssp
Q_{5,3}=\tttt{1\over M}[Q_{3,2},T_2]\scs \ssp
Q_{3,2}=\tttt{1\over M}[T_1,T_2].
$$

\noindent 
Thus
$$
Q_{7,4}=\tttt{1\over M^3}[[[T_1,T_2],T_2],T_2]
\ses
\tttt{1\over M^3}(T_1T_2^3-3T_2T_1T_2^2+3T_2^2T_1T_2-T_2^3T_1)
$$
We should mention  that it follows from the Stanley-Macdonald  Pieri rules  that to compute the action of an operator $Q_{m,n}$ we only need   its ``{\ita symbol }''  $\Xi_{m,n}$.  This is the  polynomial in $\xon$ that is obtained by replacing,
in each monomial,
 the $i^{th}$ factor $T_{a_i}$  by $x_i^{a_i}$. 
 For instance for $Q_{7,4}$  this gives
 $$
\Xi_{m,n}[x_1,x_2,x_3,x_4]\ses \tttt{1\over M}
\Big(x_1x_2^2x_3^2x_4^2
-3x_1^2x_2x_3^2x_4^2
+3x_1^2x_2^2x_3x_4^2
+x_1x_2^2x_3^2x_4^2\Big)
 $$
 In the general case, denoting by $\CS_k$ the operation of making the replacements  
 $x_i\RA x_{i+k}$ in a polynomial in the $x_i$ variables, we can  construct
 $\Xi_{m,n}[\xon]$ by the recursion
 $$
\Xi_{m,n}= 
\cases {
\tttt{1\over M}\Big(\Xi_{c,d}\CS_d\Xi_{a,b}\sms \Xi_{a,b}\CS_b\Xi_{c,d}
\Big)& if $n>1$
 and $Split(m,n)=(a,b)+(c,d)$\cr \cr
x_1^m & if $n=1$.
\cr
}
\eqno 1.17
$$
This given, as a corollary of Theorem I.3 we will  obtain that
$$
 Q_{m,n} (-1)^n=
\sum_{\mu\part n}{\TH_\mu[X;q,t]\over w_\mu}\bbu
\sum_{T\in ST(\mu)}\bbu \Xi[ w_T(1),w_T(2),\ldots ,w_T(n)]
\prod_{k=2}^n \tttt{ 1-w_T(k) qt    \over 1-qt }
\bbu \prod_{1\le h<k \le n}\bbu
f\big(\tttt{  w_T(h)\over w_T(k) }\big)
\eqno 1.18
$$

Our first task  will be  to establish the identities that permit our definition of the operators $Q_{u,v}$ when $u,v$ are not co-prime. This was carried out in [3] by viewing $\CA$ as generated by the family $\{ Q_{1,n}\}_{n\ge 0}$. In the present development, (due to  Theorem 1.1), we view 
$\CA$ as generated by the family $\{ Q_{m,1} \}_{m\ge 0}$. It turns out that an identity established in 2008
in the research that led to results in [11] 
turns out to provide the basic ingredient needed in the present development. It may be stated as follows
\sas

\noindent{\bol Theorem 1.2}

{\ita For all $m\ge 1$ the operators in the family
$$
\big\{[ Q_{b,-1},Q_{a,1}]\big\}_{\multi{a\ge 0;b\ge0 \cr a+b=m }}
\eqno 1.19
$$
all commute with $D_0$ and act identically on $\La$.}
\supereject

The proof will be given in the next section. Here we will prove what is necessary to define the operators $Q_{b,-1}$ and derive some of their properties. 

To begin, for any co-prime pair $(m,n)$ we will set
$$
Q_{m,-n}\ses Q_{m,n}^{\perp ^*}
\eqno 1.20
$$
where ``$\perp^*$'' denotes the operation of the taking the 
adjoint of a symmetric function operator with respect to the
$*$-scalar product. In other words $Q_{m,-n}$ is the unique
operator which satisfies
$$
\LL Q_{m,-n}f\scs g\RR_* =
\LL f\scs  Q_{m,n} g\RR_* 
\bigsp \hbox{(for all $f,g\in \La$)}.
\eqno 1.21
$$

To obtain an explicit formula for $Q_{b,-1}$ our starting point are the following two auxiliary identities
\sas
\noindent{\bol Lemma 1.1}

{\ita For all $f,g\in \La$  we have
$$
\LL f\sscs g\RR_* \ses \LL \phi\, \om f\sscs g\RR 
\ses \LL \om\, \phi\, f\sscs g\RR 
$$
where $\phi$ is the operator defined by the plethystic substitution}
$$
\phi \, f[X]\ses f\big[MX\big]\ses f\big[{\textstyle{  (1-t)(1-q)X}}\big].
\eqno 1.22 
$$

This is an easy consequence of the definition in I.19  
(see [6] for a proof).

\noindent{\bol Lemma  1.2}
$$
Q_{0,-1}\ses -M e_1^\perp
\eqno 1.23
$$
{\ita where $e_1^\perp $ is the adjoint of $\ul e_1$ with respect to the Hall scalar product.}

\noindent{\bol Proof}

Recall from 1.13 a) that $Q_{0,1}=-e_1$ thus by 1.22
$$
\LL Q_{0,1} f\scs g\RR_*=-\LL e_1 f\scs g\RR_*
=
- \LL \phi\om e_1 f\scs g\RR=-M \LL  e_1 \phi\om f\scs g\RR
=-M \LL    \phi\om f\scs e_1^\perp g\RR
=-M \LL      f\scs e_1^\perp g\RR_*.
$$

As a corollary we get
\sas

\noindent{\bol Proposition  1.3}

{\ita For all $m\ge 1$  we have}
$$
Q_{m,-1}\ses -M \nabla^{-m}e_1^\perp\nabla^m.
\eqno 1.24
$$
\noindent{\bol Proof}

From the particular case $n=1$ of 1.20  we have
$$
\LL Q_{m,-1}f\scs g\RR_*=\LL f \scs Q_{m,1} g\RR_*
\bigsp \hbox{(for all $f,g\in \La$)}
$$
and 1.15 gives
$$
\LL Q_{m,-1}f\scs g\RR_*=\LL f \scs \nabla^m Q_{0,1}\nabla^{-m} g\RR_*=-M \LL \nabla^{-m}e_1^\perp\nabla^m f \scs   g\RR_*.
$$
The last equality is due to the self-adjointness of $\nabla$ with respect to the $*$-scalar product.
\sas

The  following identities will also play a role in the sequel
\sas

\vbox{
\noindent{\bol Proposition  1.4}

{\ita For all $a,b\ge 1$  we have}
$$
a) \ess \ess Q_{a,1}=\tttt{1\over M}\big[Q_{a-1,1}\scs D_0\big]
\bigsp\bigsp
b) \ess \ess Q_{b,-1}=\tttt{1\over M}\big[D_0\scs Q_{b-1,-1}  \big].
\eqno 1.25
$$}
\noindent{\bol Proof}

The identity in 1.25 a) is an instance of 1.16. The identity in 1.25 b) is the $*$-adjoint of 
1.25 a) together  with the fact that taking ``$\perp^*$'' reverses order.
\sa

For  $F[X;q,t]\in \Lambda$ let us set
$$
\DA\,F[X;q,t\,]\ses \om\, F[X;1/q,1/t\,] \ses F[-\eee X;1/q,1/t\,]\ess .
\eqno 1.26
$$
It is easily seen that the operator ``$\DA$'' is an involution. It also has 
the following useful properties, proved in [6].
\sas
 
\noindent{\bol Proposition 1.5}

$$
\eqalign{
&a)\ess \ess \ess \DA\, \nabla\,  \DA  \ses \nabla ^{-1}
\cr
& b)\ess \ess \ess \DA \, D_k\,  \DA  \ses (-1)^k\, D_k^*\ess .
\cr}
\eqno 1.27
$$

\sas

\noindent{\bol 2. Proofs of  basic identities for the algebra $\CA$}
\sa

Our  goal in this section  is to prove  Theorem I.1 
and its corollary Theorem I.2.  This requires  establishing 
first Theorem 1.2. To carry all this out, we need some preliminary
observations and establish some auxiliary properties of the family of operators
$$
Q_{m,-1}= \nabla^{-m}Q_{0,-1}\nabla^m= -M \nabla ^{-m} e_1^\perp \nabla^m.
\eqno 2.1
$$
For notational convenience we need to set 
$$
 \nabla^{-m}e_1^\perp \nabla^m\ses V_m.
\eqno 2.2
$$
This given, 1.24 and and 1.25 b), namely the two identities 
$$
Q_{m,-1}=-M\nabla^{-m}e_1^\perp \nabla^m \scs
\bigsp
Q_{b,-1}=\tttt{1\over M}[D_0,Q_{b-1,-1}]
$$
combine to give us the recursion
$$
V_m\ses \tttt{1\over M}[D_0\scs V_{m-1}].
\eqno 2.3
$$

Surprisingly, a simple conjugation by the ``$\DA$'' operators 
reverses this recursion. More precisely we have 
\sas

\noindent{\bol Proposition 2.3}

$$
V_{m-1}\ses  \tttt{1\over \tM}  [D_0^*,V_m] .
\eqno 2.4
 $$
 \noindent{\bol Proof}

The definition in 2.2 and 2.3 for $m=1$ give
$$
\nabla^{-1}e_1^\perp \nabla \ses 
\tttt{1\over M}\big[D_0\scs   e_1^\perp   \big].
$$
Since we trivially have $\DA e_1^\perp \DA=e_1^\perp$ (as it is easily verified by applying both sides to any of the standard symmetric function bases), from Proposition 1.4 we derive that
$$
\nabla^{  }e_1^\perp \nabla^{-1 }\ses 
\tttt{1\over \tM}\big[D_0^*\scs  e_1^\perp  \big ].
\eqno 2.5
$$
Since   $\nabla$ and $D_0^*$  are both eigen-operators for the modified  Macdonald basis $\{\TH_\mu[X;q,t]\}_\mu$, they commute. Thus conjugating both sides of 2.5
by $\nabla^{-m}$ gives
$$
\nabla^{ -m+1 }e_1^\perp \nabla^{m-1 }\ses 
\tttt{1\over \tM}\big[D_0^*\scs \nabla^{-m} e_1^\perp \nabla^{m} \big ],
$$
which is another way of writing 2.4.
\sas

\noindent{\bol Theorem 2.1}

{\ita The operator
$$
U_m\ses \big[ \nabla ^{-m} e_1^\perp \nabla^m\scs \ul e_1\big]
\eqno 2.6
$$
commutes with both $D_0$  and $D_0^*$ and therefore it   is
an eigen-operator  of the   basis 
$\big\{ \TH_\mu [X;q,t]\big\}_\mu$.}

\noindent{\bol Proof}

We will prove by induction that we have for all $m\ge 1$
$$
\cases
{
a)\ess\ess\ess U_m=[D_1,V_{m-1}]\cr\cr
b)\ess\ess\ess [U_m,D_0^*]=0\cr
}.
\eqno 2.7
$$
Note that setting  $m=1$ in 2.6 and using 1.11 $(iii)$ and then 
1.11 $(ii)$ we get
$$
U_1\ses  \big[ \nabla ^{-1} e_1^\perp \nabla \scs \ul e_1\big]
\ses\tttt{1\over M} \big[D_{-1} \scs \ul e_1\big]
\ses D_0.
$$
This proves 2.7 b) for $m=1$. To prove 2.7  in the base case, we are left to show that
$$
D_0\ses [D_1,e_1^\perp],
$$
but this is precisely  1.11$(v)$. We can thus  inductively assume 2.7 true up to $m-1$.


Now by definition we have
$$
U_{m}\ses [V_{m},e_1]
$$
and 2.3 gives
$$
\eqalign{
U_{m}
&\ses  \tttt{1\over M} \big[[D_0,V_{m-1}],e_1\big]
\cr \hbox{(by Jacobi)}
&\ses  - \tttt{1\over M} \big[[ V_{m-1},e_1],D_0\big]\sms  \big[[e_1,D_0],V_{m-1}\big] 
\cr \hbox{(by 1.11  $(ii)$)} 
&\ses  - \tttt{1\over M} \big[U_{m-1},D_0\big]\sps  \big[D_1,V_{m-1}\big] 
\cr \hbox{(by 2.7 b) for $m-1$)} 
&\ses   \big[D_1,V_{m-1}\big] .
\cr}
$$
This proves 2.7 a)  for $m$. To show 2.7 b)  for $m$ we use this and get
$$
\eqalign{
[U_m,D_0^*]
&\ses \big [ [D_1,V_{m-1}],D_0^*,\big]
\cr \hbox{(by Jacobi)}
&\ses - \big[[V_{m-1},D_0^*],D_1\big]
\sms    \big[[D_0^*,D_1],V_{m-1}\big].
\cr}
\eqno 2.8
$$
For the first  term Proposition 2.3 gives (using induction) 
$$
- \big[[V_{m-1},D_0^*],D_1\big]
\ses
 \tM  \big[V_{m-2},D_1\big]\ses
 -\tM U_{m-1}.
\eqno 2.9
 $$
For the second  term  in 2.8 we have, using 1.11 $(iii)$, 
$$
\eqalign{
 [D_0^*,D_1] 
 &=\big[\nabla e_1\nabla^{-1}\scs D_0^*\big]\cr
  (\hbox{$D_0^*$ and $\nabla$ commute})
  &=\nabla \big[ e_1\scs D_0^*\big]\nabla^{-1}\cr
   (\hbox{by 1.11 $(ii)^*$})
  &=\tM \nabla D_1^*\nabla^{-1}\cr  
    (\hbox{by 1.11 $(iii)^*$})
  &=\tM \ul e_1.\cr  
\cr
}
$$
\vskip -.28in
\noindent
Thus
$$
\sms    \big[[D_0^*,D_1],V_{m-1}\big]
\ses -\tM \big[e_1,V_{m-1}\big]\ses \tM U_{m-1}
$$
and this together with 2.9  reduces 2.8 to
$$
[U_m,D_0^*]\ses -\tM U_{m-1}+\tM U_{m-1}\ses 0
$$
completing the induction.
\sas

We are now finally in a position to prove the following sharpening of Theorem 1.2.
\sas
\noindent{\bol Theorem 2.2}

{\ita For all $m\ge 1$ the operators in the family
$$
\big\{[ Q_{b,-1},Q_{a,1}]\big\}_{\multi{a\ge 0;b\ge0 \cr a+b=m }}
\eqno 2.10
$$
 act identically on $\La$ and we can set
 $$
 Q_{m,0}\ses  \tttt{1\over M}\big[Q_{m,-1},Q_{0,1} \big]\ses 
 \big[  \nabla ^{-m} e_1^\perp \nabla^m\scs
\ul e_1\big].
 \eqno 2.11
 $$
 In particular we see that the collection $\{Q_{m,0}\}_{m\ge0}$
 is a commuting family of eigen-operators for the modified Macdonald basis.}

\noindent{\bol Proof}

Notice first that the second equality in 2.11 follows from  2.1 and 1.13 a). Moreover, since Theorem 2.1 assures that the right hand side of 2.11 commutes with $D_0$
and $D_0^*$, the last assertion follows from the fact that the eigenvalues of $D_0$ (or $D_0^*$)    are all distinct. Thus the only thing that remains to prove is the equality of these operators. 
But this is now easily seen if we write the operators in 2.10
using 1.15  and 1.24, that is

$$
[ Q_{b,-1},Q_{a,1}]
\ses  M
\big[ \nabla ^{-b} e_1^\perp \nabla^b,  \nabla^a e_1\nabla^{-a}\big].
$$
In fact,  we know from Theorem 2.1 that this operator, for $b=m$ and $a=0$, commutes with $\nabla$ thus
$$
[ Q_{m,-1},Q_{0,1}] 
=M\nabla \big[ \nabla ^{-m } e_1^\perp \nabla^{m },   e_1 \big]\nabla^{-1}
=M\big[ \nabla ^{-m+1} e_1^\perp \nabla^{m-1},  \nabla  e_1\nabla^{-1}\big]
=[ Q_{m-1,-1},Q_{1,1}].
$$
Proceeding by descent induction on $b\in[1,m]$, assume that for  $a=m-b$ we have
$$
[ Q_{m,-1},Q_{0,1}] =[ Q_{b,-1},Q_{a,1}].
$$
Conjugating by $\nabla$ we similarly obtain
$$
[ Q_{m,-1},Q_{0,1}] =\nabla [ Q_{b,-1},Q_{a,1}]\nabla^{-1}=  [ Q_{b-1,-1},Q_{a+1,1}].
$$
This completes the induction and our proof.
\sa

Our next and final task in this section is the identification of the eigenvalues of the operators 
$$
Q_{m,0}=\big[ \nabla ^{-m} e_1^\perp \nabla^m\scs \ul e_1\big].
$$
 This task, as well as the
developments in the next section are heavily dependent on
the Pieri rules for $\ul e_1$ and $e_1^\perp$ and their summation formulas. These remarkable identities may be stated as follows.

\sas \sas

\noindent{\bol Theorem 2.3}

{\ita For  $\mu\part n$ and $\nu\part n-1$ we have
$$
a)\ess\ess\ess e_1\TH_\nu\ses \sum_{\mu\LA \nu }d_{\mu,\nu}\TH_\mu
\scs \bigsp
b)\ess\ess\ess e_1^\perp\TH_\mu\ses \sum_{\nu\RA \mu }c_{\mu,\nu}\TH_\nu
\eqno  2.12
$$
where ``$\nu\RA \mu$''   means that $\mu$ is obtained 
by adding a corner square to $\mu$. Moreover we also have the basic relation
$$
d_{\mu,\nu}\ses M  {w_\nu\over w_\mu}c_{\mu,\nu}.
\eqno 2.13
$$
}

\noindent{\bol Theorem 2.4} 
$$
\eqalign{
a)\ess\ess\ess \sum_{\nu\RA\mu}c_{\mu\nu}(q,t)\, (T_\mu/T_\nu)^k
&\ses 
\cases
{{tq\over M}\ssp h_{k+1}\big[D_\mu(q,t)/tq\big] & if $k\geq 1\ess  $\cr\cr
B_\mu(q,t) & if $k=0\ess  $
\cr
}
\cr\cr
b)\ess\ess\ess  \sum_{\mu\LAr\nu}d_{\mu\nu}(q,t)\, (T_\mu/T_\nu)^k
&\ses 
\cases
{(-1)^{k-1}\ssp e_{k-1}\big[D_\nu(q,t) \big] & if $k\geq 1\ess  $\cr\cr
1 & if $k=0\ess   $
\cr
}
\cr
}.
\eqno 2.14
$$

\noindent{\bol Remark 2.1}

We must mention that  2.12 b) and 2.13 a) were 
first proved in [10], directly from the original Stanley-Macdonald Pieri rules. On the other hand  for 2.12 a),  and 2.13 b),  a direct proof was    never published, 
although extensively used in several publications. These proofs will be included in the next section for sake of completeness. An indirect derivation of 2.14 a) and b) was also given
in [6].
\sas

Our goal, which  is  Theorem I.1,   may be simply restated as follows
\sas

\noindent{\bol Theorem 2.5 }

{\ita For all $k\ge 1 $ and partitions $\mu$ we have}
$$
\big [ \nabla^{-k}e_1^\perp \nabla^k\scs  e_1\big]\TH_\mu[X;q,t] \ses
 \tttt{qt \over qt-1}  h_k \big[ D_\mu(q,t) (\tttt{1\over qt}-1) \big]\TH_\mu[X;q,t].
\eqno 2.15
$$

\noindent{\bol  Proof}

Notice first  that if we know already that a certain symmetric function operator $\BZ$ is an eigen-operator for the basis $\{\TH_\mu[X;q,t]\}_\mu$, 
then  given the expansion
$$
e_n^*[X]=e_n[\tttt{X\over M}]\ses \sum_{\mu\part n}{\TH_\mu[X;q,t]\over w_\mu}
$$
 we can simply identify its eigenvalues ${\bf z}_\mu$ from the formula
$$
\BZ \, e_n[\tttt{X\over M}]\ses \sum_{\mu\part n}{\TH_\mu[X;q,t]\over w_\mu}
\,\, {\bf z}_\mu.
$$
This given, let us  start with
$$
\eqalign{
\nabla^{ k}e_1e_n^*
&\ses \sum_{\mu\part n }{1\over w_\mu}
\sum_{\ggg\LA \mu}d_{\ggg\mu}\TH_\ggg T_\ggg^k
\cr
&\ses \sum_{\ggg\part n+1 }{\TH_\ggg T_\ggg^k\over w_\ggg}
\sum_{\mu\RA \ggg}\tttt{w_\ggg\over w_\mu}d_{\ggg\mu} 
\cr
\hbox{(by 2.13)}
&\ses M\sum_{\ggg\part n+1 }{\TH_\ggg T_\ggg^k\over w_\ggg}
\sum_{\mu\RA \ggg} c_{\ggg\mu} 
\cr
\hbox{(by 2.14 a))}
&\ses M\sum_{\ggg\part n+1 }{\TH_\ggg T_\ggg^k\over w_\ggg}
B_{\ggg}(q,t).
\cr}
$$
Thus
$$
\eqalign{
\nabla^{-k}e_1^\perp \nabla^{ k}e_1e_n^*
&\ses 
M\sum_{\ggg\part n+1 }{\nabla^{-k}e_1^\perp\TH_\ggg T_\ggg^k\over w_\ggg}
B_{\ggg}(q,t) 
\cr
&\ses
 M\sum_{\ggg\part n+1 }{  T_\ggg^k\over w_\ggg}
B_{\ggg}(q,t)
\sum_{\ddd\RA\ggg}c_{\ggg\ddd}\TH_\ddd T_\ddd^{-k}
\cr
&\ses
 M\sum_{\ddd\part n  }{\TH_\ddd\over w_\ddd}
\sum_{\ggg\LA\ddd}\tttt{w_\ddd\over w_\ggg} 
\, c_{\ggg\ddd}\, \tttt{T_\ggg^k\over T_\ddd^{ k}} 
\ssp B_{\ggg}(q,t)
\ses   \sum_{\ddd\part n  }{\TH_\ddd\over w_\ddd}
\sum_{\ggg\LA\ddd} 
\, d_{\ggg\ddd}\, \tttt{T_\ggg^k\over T_\ddd^{ k}} 
\ssp B_{\ggg}(q,t).
\cr}
\eqno  2.16
$$
But the relation $B_\ggg=B_\ddd+{T_\ggg\over T_\ddd}$ gives
$$
\sum_{\ggg\LA\ddd} 
\, d_{\ggg\ddd}\, \tttt{T_\ggg^k\over T_\ddd^{ k}} 
\ssp B_{\ggg} 
\ses 
B_\ddd\sum_{\ggg\LA\ddd} 
\, d_{\ggg\ddd}\, \tttt{T_\ggg^k\over T_\ddd^{ k}} 
\ssp\sps \sum_{\ggg\LA\ddd} 
\, d_{\ggg\ddd}\, \tttt{T_\ggg^{k+1}\over T_\ddd^{ k+1}} 
$$
and 2.14  b) delivers
$$
\sum_{\ggg\LA\ddd} 
\, d_{\ggg\ddd}\, \tttt{T_\ggg^k\over T_\ddd^{ k}} 
\ssp B_{\ggg} 
\ses  (-1)^{k-1}\ssp B_\ddd\,  e_{k-1}\big[D_\ddd(q,t) \big]
\sps   (-1)^{k } \,  e_{k }\big[D_\ddd(q,t) \big].
$$
Using this in 2.16 we obtain
$$
\nabla^{-k}e_1^\perp \nabla^{ k}e_1e_n^*
\ses   (-1)^{k-1}\ssp \sum_{\ddd\part n  }{\TH_\ddd\over w_\ddd}
 B_\ddd\,  e_{k-1}\big[D_\ddd(q,t) \big]
\sps   (-1)^{k }\sum_{\ddd\part n  }{\TH_\ddd\over w_\ddd} \,  e_{k }\big[D_\ddd(q,t) \big].
\eqno  2.17
$$
Next we start with
$$
\eqalign{
\nabla^{-k}e_1^\perp \nabla^{ k}e_n^*
&= \sum_{\mu\part n}
{T_\mu^{k}\over w_\mu}\sum_{\nu\RA \mu}   c_{\mu\nu}{T_\nu^{-k}}\TH_\nu 
\cr
&=  \sum_{\nu\part n-1}{\TH_\nu \over w_\nu}
\sum_{\nu\RA \mu}   c_{\mu\nu}\tttt{w_\nu \over w_\mu}\tttt{T_\mu^{k}\over T_\nu^k}
=\tttt {1\over M}\sum_{\nu\part n-1}{\TH_\nu \over w_\nu}
\sum_{\mu\LA \nu}   d_{\mu\nu} \tttt{T_\mu^{k}\over T_\nu^k}
\cr}
$$
and 2.14 b)  gives
$$
\nabla^{-k}e_1^\perp \nabla^{ k}e_n^*
\ses 
 \tttt{(-1)^{k-1}\over M}\sum_{\nu\part n-1}{\TH_\nu \over w_\nu}
\ssp e_{k-1}\big[D_\nu(q,t) \big].
$$
Multiplying  on both sides by $e_1$ yields 
$$
\eqalign{
e_1\nabla^{-k}e_1^\perp \nabla^{ k}e_n^*
&= 
 \tttt{(-1)^{k-1}\over M}\sum_{\nu\part n-1}{1\over w_\nu}
\ssp e_{k-1}\big[D_\nu(q,t) \big]
\sum_{\ddd\LA \nu }d_{\ddd\nu }\TH_\ddd
\cr
&= 
\tttt {(-1)^{k-1}\over M}\sum_{\ddd\part n }{\TH_\ddd \over w_\ddd}
\sum_{\nu\RA \ddd }\tttt{ w_\ddd\over w_\nu}\, d_{\ddd\nu }
\ssp e_{k-1}\big[D_\nu(q,t) \big]
= 
\tttt {(-1)^{k-1} }\sum_{\ddd\part n }{\TH_\ddd \over w_\ddd}
\sum_{\nu\RA \ddd } \, c_{\ddd\nu }
\ssp e_{k-1}\big[D_\nu(q,t) \big].
\cr}
\eqno 2.18
$$
To deal with the sum
$$
 \sum_{\nu\RA \ddd } \, c_{\ddd\nu }
\ssp e_{k-1}\big[D_\nu(q,t) \big]
$$
the identity $D_\nu =D_\ddd-M {T_\ddd\over T_\nu} $ and the
 addition formula for the ${e_k}'s$ gives 
$$
\eqalign{
e_{k-1}[D_\nu]
&= e_{k-1}[D_\ddd]\sps
\sum_{r=1}^{k-1} e_{k-1-r}[D_\ddd]e_r[-M](\tttt{T_\ddd\over T_\nu})^r
=e_{k-1}[D_\ddd]\sps
\sum_{r=1}^{k-1}(-1)^re_{k-1-r}[D_\ddd]h_r[M](\tttt{T_\ddd\over T_\nu})^r.
\cr
}
$$
Thus
$$
\eqalign{
 \sum_{\nu\RA \ddd } \, c_{\ddd\nu }
 \ssp e_{k-1}\big[D_\nu(q,t) \big]
 &\ses
\sum_{\nu\RA \ddd } \, c_{\ddd\nu } e_{k-1}[D_\ddd]\sps
\sum_{r=1}^{k-1}(-1)^re_{k-1-r}[D_\ddd]h_r[M]
\sum_{\nu\RA \ddd } \, c_{\ddd\nu }(\tttt{T_\ddd\over T_\nu})^r
\cr
}
$$
and 2.14 b) gives
$$
\eqalign{
\sum_{\nu\RA \ddd } \, c_{\ddd\nu }  \ssp e_{k-1}\big[D_\nu(q,t) \big]
 &\ses
  e_{k-1}[D_\ddd] B_\ddd \sps
\sum_{r=1}^{k-1}(-1)^re_{k-1-r}[D_\ddd]h_r[M]
\tttt{qt\over M}h_{r+1}[D_\ddd/qt]
\cr
\hbox{(by 1.10)}
&=
  e_{k-1}[D_\ddd] B_\ddd +
\tttt{qt\over (1-qt)}\sum_{r=0}^{k-1}(-1)^re_{k-1-r}[D_\ddd](1-(qt)^r) 
h_{r+1}[D_\ddd/qt]
\cr
&=
  e_{k-1}[D_\ddd] B_\ddd -
\tttt{qt\over (1-qt)}\sum_{r=1}^{k }(-1)^{r }e_{k -r}[D_\ddd](1-(qt)^{r-1}) 
h_{r}[D_\ddd/qt]
\cr
& =  e_{k-1}[D_\ddd] B_\ddd -
\tttt{qt\over (1-qt)}\sum_{r=1}^{k } e_{k -r}[D_\ddd]
e_{r}[-D_\ddd/qt]
 + 
 \tttt{1\over (1-qt)}\sum_{r=1}^{k } e_{k -r}[D_\ddd]   
e_{r}[-D_\ddd ]
\cr
&=
  e_{k-1}[D_\ddd] B_\ddd 
\sms
\tttt{qt\over (1-qt)}  e_{k  }[D_\ddd(1-1/qt)] 
+\tttt{qt\over (1-qt)}  e_{k  }[D_\ddd ] -\tttt{1\over (1-qt)}  e_{k  }[D_\ddd ] 
\cr
\cr
&=
  e_{k-1}[D_\ddd] B_\ddd 
\sms
\tttt{qt\over (1-qt)}  e_{k  }[D_\ddd(1-1/qt)] 
\sms   \, e_{k }[D_\ddd].
\cr
 }
 $$
Using this in 2.18  we finally obtain
$$
e_1\nabla^{-k}e_1^\perp \nabla^{ k}e_n^*
\ses 
 {(-1)^{k-1} }\sum_{\ddd\part n }{\TH_\ddd \over w_\ddd}
 e_{k-1}[D_\ddd] B_\ddd 
\ssp \sps  { }\tttt{(-1)^{k}qt\over (1-qt)}\sum_{\ddd\part n }{\TH_\ddd \over w_\ddd}
\  e_{k  }[D_\ddd(1-1/qt)] \sps {(-1)^{k} }\sum_{\ddd\part n }{\TH_\ddd \over w_\ddd}
\, e_{k }[D_\ddd].
$$
Recalling that the identity in  2.17 is
$$
\nabla^{-k}e_1^\perp \nabla^{ k}e_1e_n^*
\ses   (-1)^{k-1}\ssp \sum_{\ddd\part n  }{\TH_\ddd\over w_\ddd}
 B_\ddd\,  e_{k-1}\big[D_\ddd \big]
\sps   (-1)^{k }\sum_{\ddd\part n  }{\TH_\ddd\over w_\ddd} \,  e_{k }\big[D_\ddd  \big],
$$
subtracting   this from the former 
yields the final identity
$$
 \eqalign{
 [ \nabla^{-k}e_1^\perp \nabla^k\scs  e_1]\,\,  e_n^*
 &\ses 
 \tttt{(-1)^{k}qt\over  qt-1 }\sum_{\ddd\part n }{\TH_\ddd \over w_\ddd}
\  e_{k  }[D_\ddd(1-1/qt)] 
\ses 
 \tttt{ qt\over  qt-1}\sum_{\ddd\part n }{\TH_\ddd \over w_\ddd}
\  h_{k  }[D_\ddd(\tttt{1\over qt}-1)] 
\cr}
$$
completing our proof of 2.15.
\sa

\noindent{\bol Remark 2.2}

The identity of Proposition 2.3, namely
$$
[D_0^*,V_k]\ses \tM \, V_{k-1},
\eqno 2.19
$$
is quite remarkable when viewed geometrically. In fact from 2.2 and 2.1 we derive that
$$
V_m\ses  \nabla^{-m} e^\perp  \nabla^{m}\ses
-\tttt{1\over M} Q_{m,-1}
$$

\vskip -.2in
\noindent
and the definition in  I.11 gives
$$
 D_0^*\ses -\tttt{1\over qt}\, Q_{-1,0} 
$$

\vskip -.2in
\noindent
Thus 2.19  is none other than
$$
 -\tttt{1\over qt}[Q_{-1,0},  Q_{m,-1}]\ses  \tttt{M\over qt}  Q_{m-1,1}
$$

\vskip -.3in
\noindent
or better
\vskip -.2in
$$
Q_{m-1,1}\ses \tttt{1\over M} [Q_{m,-1}, Q_{-1,0}  ],
\eqno 2.20
$$
which may be viewed as  the identity resulting from the following splitting of the vector $(m-1,-1)$.
\vskip -.2in
$$
\vcenter{\hbox{\includegraphics[width=2in]{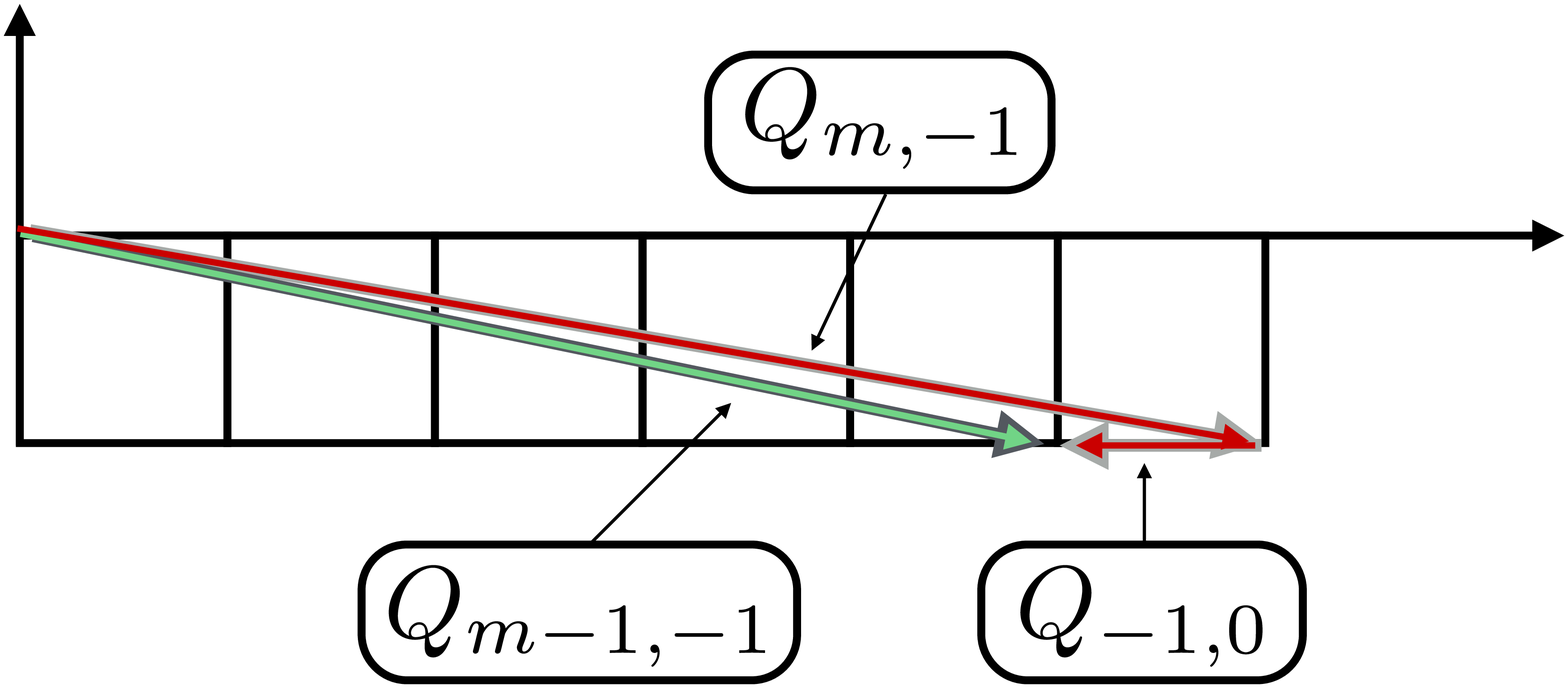}}}
$$

\noindent{\bol 3. Pieri Rules and  Standard Tableaux expansions}
\sas

 Explicit formulas for the coefficients $d_{\mu,\nu}$ 
  and $c_{\mu,\nu}$ in 2.12
were first obtained in [10] from   Macdonald's formula for the multiplication by $e_1$
of his original $P_\mu[X;q,t]$ basis. More precisely the latter formula was used in [10]  to obtain 
the identity

\vskip -.2in
$$
d_{\mu\nu}(q,t)\ses
\prod_{s\in R_{\mu        \nu}} 
{ 
q^{a_\nu(s)}-t^{l_\nu(s)+1}
\over 
q^{a_\mu(s)}-t^{l_\mu(s)+1}
}
 \prod_{s\in C_{\mu\nu}} 
{ 
t^{l_\nu(s)}-q^{a_\nu(s)+1}
\over 
t^{l_\mu(s)}-q^{a_\mu(s)+1}
}
\eqno 3.1
$$

\vskip -.2in
\noindent
where $R_{\mu\nu}$ and  $ C_{\mu\nu}$  denote the collections of cells of $\nu$ that 
are respectively in the row and the column of the cell $\mu/\nu$.
This done,  an easy use of the  orthogonality relations in I.18 gave that
$$
c_{\mu\nu}(q,t)\ses \tttt{1\over M}\tttt{w_\mu\over w_\nu} d_{\mu\nu}(q,t).
\eqno 3.2
$$
Finally, this relation  combined with 3.1, after many  cancellations, yielded  the identity
$$
c_{\mu\nu}(q,t)= \prod_{s\in R_{\mu\nu}} 
{ 
t^{l_\mu(s)}-q^{a_\mu(s)+1}
\over 
t^{l_\nu(s)}-q^{a_\nu(s)+1}
}
\prod_{s\in C_{\mu\nu}} 
{ 
q^{a_\mu(s)}-t^{l_\mu(s)+1}
\over 
q^{a_\nu(s)}-t^{l_\nu(s)+1}
}.
\eqno 3.3
$$

\vskip -.2in
It is not difficult to see that, in both 3.1 and 3.2, there are
 still many  cancellations remaining.
This observation led to   extremely useful   formulas for these Pieri coefficients.
To state them we need notation.
For a partition $\mu$ with $\l$ removable corners, let 
$x_0,x_1,x_2,\ldots,x_\l$ and $u_0,u_1,u_2,\ldots,u_\l$
denote  

\hsize 4.7in 

\noindent
the weights \hskip -.05in
\footnote{$\ ^{(\dag)}$}
{$weight(c)= t^{l_\mu'(c)}q^{a_\mu'(c)}$} 
  of the cells of $\mu$ as illustrated  in the adjacent  diagram  for $\l=5$. 
It is shown in    [10]    that  these cancellations   reduce  3.3   to
$$
c_{\mu,\nu^{(k)}}\ses {x_k\over (1-\tttt{1\over t})(1-\tttt{1\over q})}
{\prod_{i=0}^\l\big(1-\tttt{u_i\over x_k}\big)\over  \prod_{i=1;i\neq k}^\l\big(1-\tttt{x_i\over x_k}\big)},
\eqno 3.4
$$

\hsize 6.5in
\vskip -.2in
\noindent
where $\nu^{(k)}$ is the partition obtained  by removing from $\mu$ the corner with weight $x_k= T_\mu/T_{\nu^{(k)}}$.

\vskip -1.4in
 \hfill $
\vcenter{\hbox{\includegraphics[width=1.6 in]{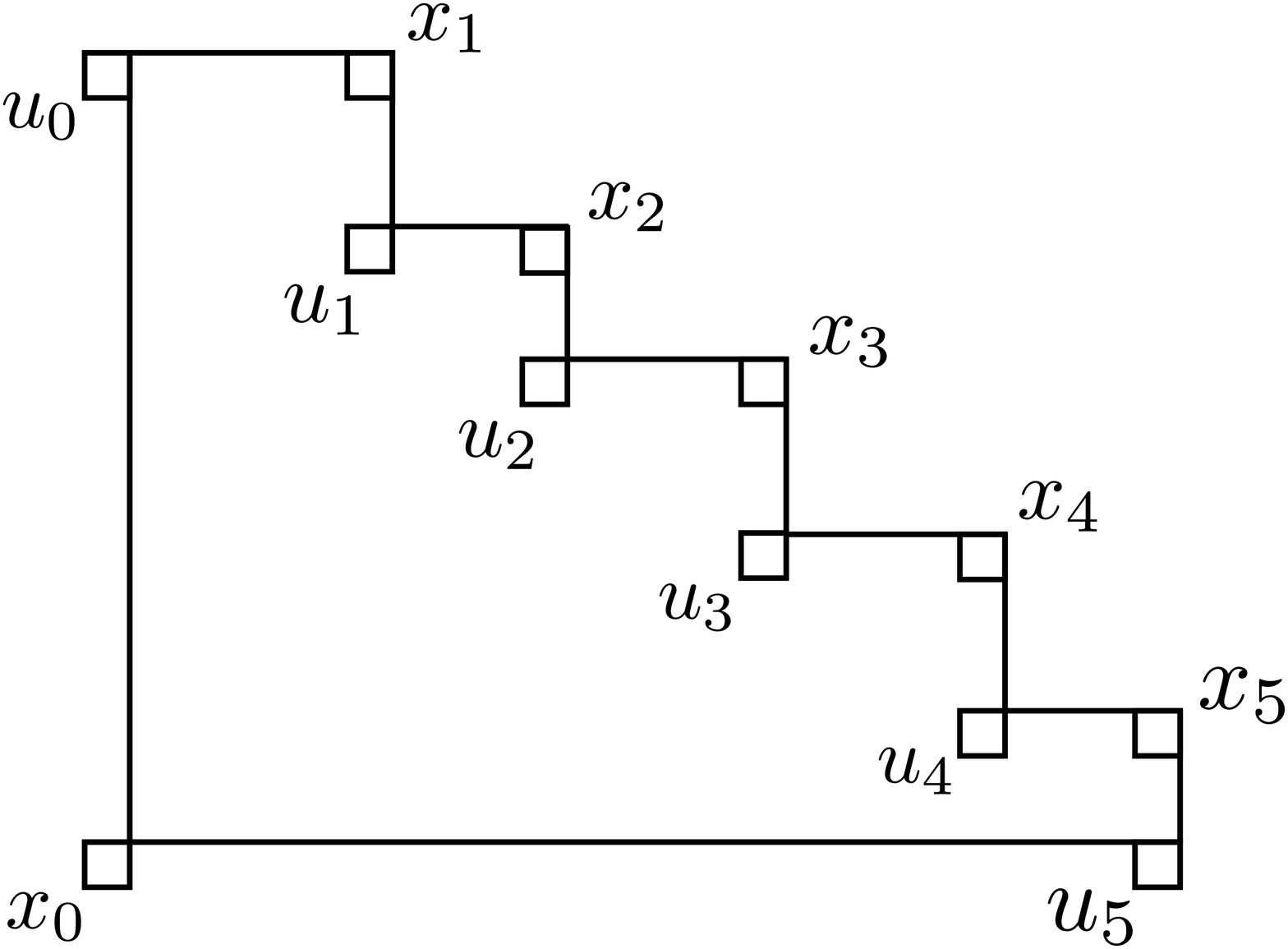}}}
 $
\supereject
\vskip .25in
 Moreover, it is shown  in [10] that massive cancellations also yield the truly remarkable identity
$$
(1-\tttt{1\over t})(1-\tttt{1\over q})B_\mu(q,t)
\ses x_0+x_1+\cdots +x_\l\sms u_0-u_1-\cdots -u_\l.
\eqno 3.5
$$

Somewhat later the last named author 
derived the companion formula for the  Pieri coefficients $d_{\mu,\nu}$. This formula, 
which  can be obtained by carrying out the appropriate cancellations in 3.1, can be written as
$$
d_{\mu^{(k)},\nu}\ses {1\over q\, t\, u_k} 
{\prod_{i=1}^\l\big(1-\tttt{x_i\over u_k}\big)\over  \prod_{i=0;i\neq k}^\l\big(1-\tttt{u_i\over u_k}\big)}
\eqno 3.6
$$
provided the shape in   3.4 is now interpreted as the Ferrers' diagram of the partition $\nu$
and $\mu^{(k)}$ is now the partition that is obtained from $\nu$ by adding the (addable)  corner square 
that is  the NE shift by one cell unit of the cell labelled by $u_k$ in 3.4. In particular this gives that
$$
{T_{\mu^{(k)}}/ T_\nu}\ses t\,q \, u_k.
\eqno 3.7
$$
Since the proof of 3.6 was never published, it will be good to include it here
and at the same time illustrate the process that yields 
3.6 from 3.1.

To begin let us start from the figure in 3.4 interpreted as the Ferrers' diagram of $\nu$ but shift all the labelled cells NE by one cell unit,  placing  a bar on each of their labels. That is 
we are setting $\overline x_i=qtx_i$ and $\overline u_i=qt u_i$. With these conventions,
multiplying both sides of the identity in  3.5 by $qt$ and replacing $\mu$ by $\nu$ we obtain 
 $$
MB_\nu(q,t)\ses  \ol x_0+\ol x_1+\cdots +\ol x_\l\sms \ol u_0-\ol u_1-\cdots -\ol u_\l
\eqno 3.8
$$ 

In the following display we have on the right the labelled Ferrers' diagram of $\nu$
and on the left   we have 3.1 atop the identity we will  derive from it. 
$$ 
{
\displaystyle{d_{\mu\nu}(q,t)\ses
\prod_{s\in R_{\mu\nu}} 
{ 
q^{a_\nu(s)}-t^{l_\nu(s)+1}
\over 
q^{a_\mu(s)}-t^{l_\mu(s)+1}
}
 \prod_{s\in C_{\mu\nu}} 
{ 
t^{l_\nu(s)}-q^{a_\nu(s)+1}
\over 
t^{l_\mu(s)}-q^{a_\mu(s)+1}
}}
\atop
d_{\mu^{(k)},\nu}\ses {1\over \ol u_k} 
{\prod_{i=1}^\l\big(1-\tttt{\ol x_i\over \ol u_k}\big)\over  \prod_{i=0;i\neq k}^\l\big(1-\tttt{\ol  u_i\over \ol u_k}\big)}
}
\bigsp \vcenter{\hbox{\includegraphics[width=1.7 in]{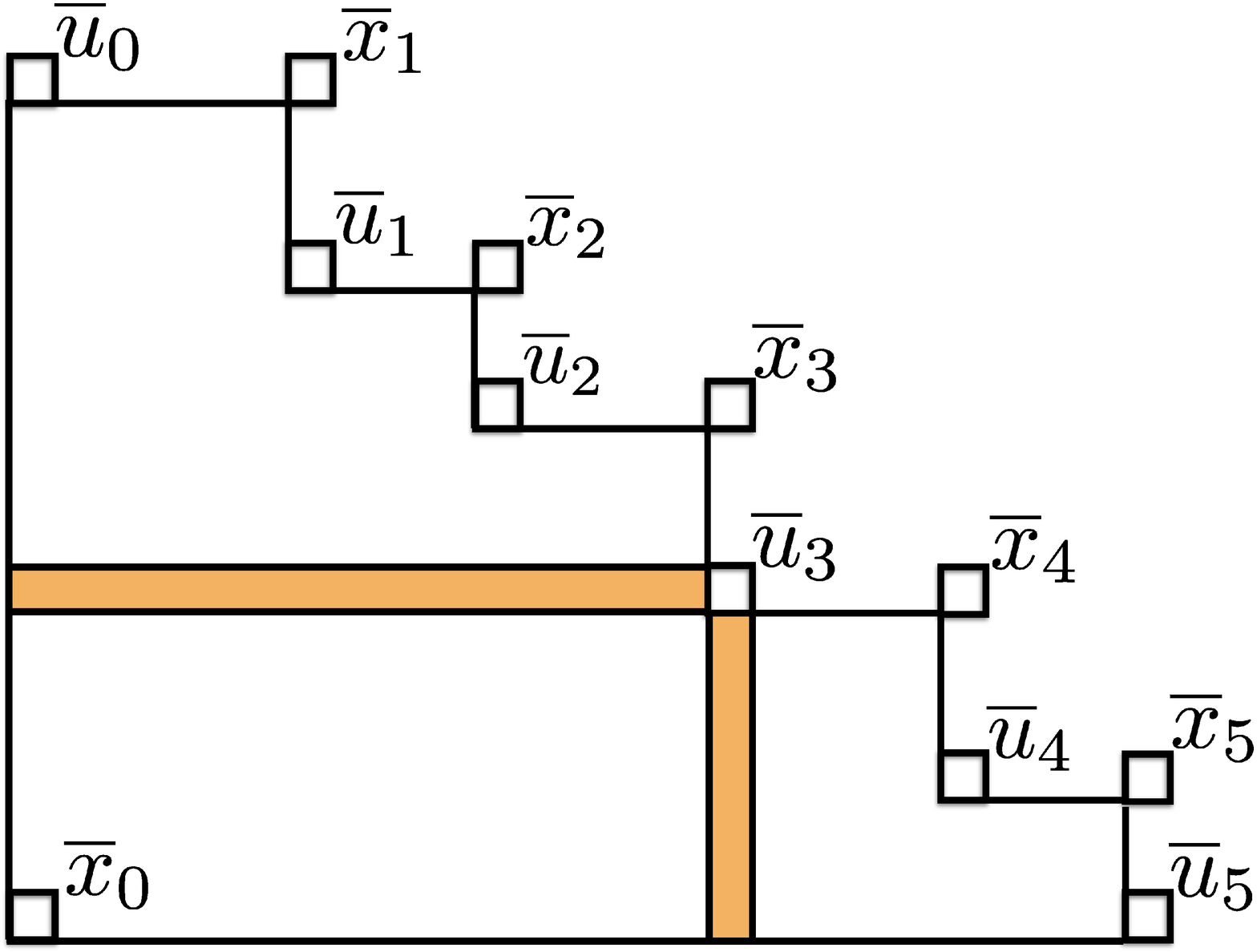}}}.
\eqno 3.9
$$
Here, for convenience, we set
$$
\ol x_i=t^{\aaa_i} q^{\bbb_i}\ess\ess\& \ess\ess
\ol u_i=t^{\aaa_{i+1}} q^{\bbb_i}
\ess\ess\ess\ess \ess(\hbox{for $0\le i\le \l $, with $\bbb_0=\aaa_{\l+1}=0$}).
\eqno 3.10
$$

To get across the cancellations that occur in S.21 along the row $R_{\mu\nu}$, we need only focus our attention on the following four figures. To help visualize where these figures are located in the diagram of $\nu$ we have depicted in  3.9 the row $R _{\mu^{(k)},\nu}$ when the cell added to  $\nu$ to obtain $\mu$ is the one whose weight is $\overline u_3$ 
$$
\vcenter{\hbox{\includegraphics[width=3.5 in]{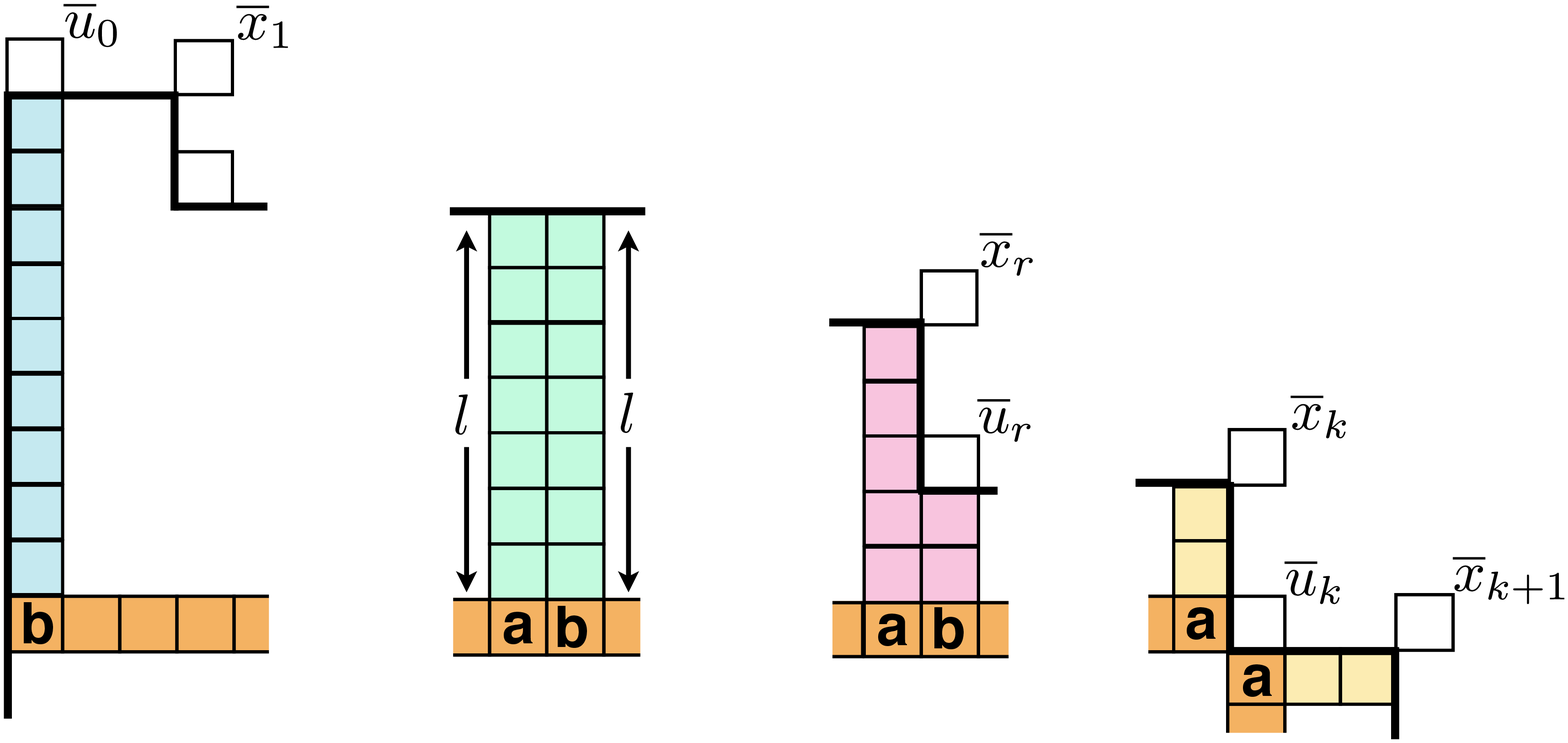}}}
$$
To begin note that for two cells $a,b\in R _{\mu^{(k)},\nu}$,  both of  whose legs are  $\l$ (as in the second figure above),   we have
$$
q^{a_\nu(a)}-t^{l_\nu(a)+1}\ses q^{(\bbb_k-1)-a}-t^{l+1}
\ses  q^{\bbb_k-b}-t^{l+1}\ses q^{a_\mu(b)}-t^{l_\mu(b)+1}.
$$ 
Thus these two factors cancel each other in  $d_{\mu^{(k)},\nu}$.

On the other hand for two cells  $a,b\in R _{\mu^{(k)},\nu}$, with $a$  preceding  $b$  and $b$ below the cell labelled $\ol u_r$ (as in the third figure above), we have
$$
q^{a_\nu(a)}-t^{l_\nu(a)+1}= q^{ \bbb_k-\bbb_r}-t^{\aaa_r-\aaa_{k+1}}
\ess\ess\ess\ess\ess \hbox{and}\ess\ess\ess\ess\ess
q^{a_\mu(b)}-t^{l_\mu(b)+1}= q^{ \bbb_k-\bbb_r}-t^{\aaa_{r+1}-\aaa_{k+1}}.
$$
Thus these two cells of $R _{\mu^{(k)},\nu}$ contribute the following ratio to $d_{\mu^{(k)},\nu}$
$$
{q^{ \bbb_k-\bbb_r}-t^{\aaa_r-\aaa_{k+1}}
\over
q^{ \bbb_k-\bbb_r}-t^{\aaa_{r+1}-\aaa_{k+1}}
}
\ses {
q^{ \bbb_k}t^{\aaa_{k+1}} -t^{\aaa_r }q^{\bbb_r}
\over
q^{ \bbb_k}t^{\aaa_{k+1}}-t^{\aaa_{r+1}}q^{\bbb_r}}\ses 
{
\ol u_k -\ol x_r
\over 
\ol u_k -\ol u_r}.\eqno 3.11
$$

These cases take care of all the factors contributed by  $ R _{\mu^{(k)},\nu}$ except for the first
factor in the denominator
 and the last factor in the numerator, respectively  contributed
by the cell $b$ in the first figure above  and the first cell $a$ in the last figure above.
These two factors yield the ratio
$$
{q^{a_\nu(a)}-t^{l_\nu(a)+1}
\over
q^{a_\mu(b)}-t^{l_\mu(b)+1}}
=
{q^{ \bbb_k-\bbb_k}-t^{\aaa_k-\aaa_{k+1}}
\over
q^{ \bbb_k-\bbb_0}-t^{\aaa_{1}-\aaa_{k+1}}
}
=
{q^{\bbb_0}
\over
q^{\bbb_k}
}{
q^{ \bbb_k}t^{\aaa_{k+1}} -t^{\aaa_k }q^{\bbb_k}
\over
q^{ \bbb_k}t^{\aaa_{k+1}}-t^{\aaa_{1}}q^{\bbb_0}}
=
{1\over q^{\bbb_k}}
{
\ol u_k -\ol x_k
\over 
\ol u_k -\ol u_0}. \eqno 3.12
$$

\hsize=4.5in
We thus obtain that
$$
\prod_{s\in R_{\mu\nu}} 
{ 
q^{a_\nu(s)}-t^{l_\nu(s)+1}
\over 
q^{a_\mu(s)}-t^{l_\mu(s)+1}
}
\ses
{1\over q^{\bbb_k}}
{\prod_{i=1}^{k}(\ol u_k-\ol x_i)
\over
{\prod_{i=0}^{k-1}(\ol u_k-\ol u_i)}
}.
\eqno  3.13
$$
To compute the contribution of the column $ C _{\mu^{(k)},\nu}$ to  $d_{\mu^{(k)},\nu}$
we start as in the previous case to note that for two adjacent cells $a,b$  with the same
arm (as indicated in the third figure  on the right), the ratio
$$
{t^{l_\nu(a)}- q^{a_\nu(a)+1}
\over
t^{l_\mu(b)}-q^{a_\mu(b)+1}
}
$$
contributes nothing, since $l_\nu(a)=l_\mu(b)$ in this case.

For the two cells $a,b\in C _{\mu^{(k)},\nu}$
as indicated in the second figure 

\hsize6.5 in

\noindent
on the right, reasoning as we did above we see that their contribution to  $d_{\mu^{(k)},\nu}$ is the ratio

$$
{t^{l_\nu(a)}- q^{a_\nu(a)+1}
\over
t^{l_\mu(b)}-q^{a_\mu(b)+1}
}
=
{t^{\aaa_{k+1}-\aaa_{r+1}} -q^{\bbb_{r+1}-\bbb_k}
\over
t^{\aaa_{k+1}-\aaa_{r+1}} -q^{\bbb_{r}-\bbb_k}
}
={
\ol u_k-\ol x_{r+1}
\over
\ol u_k-\ol u_r
}
$$
Finally the bottom cell $b\in C _{\mu^{(k)},\nu}$  and its top cell $a$, as indicated in the bottom
and top figure in the above display,   contribute the  ratio
$$
{t^{l_\nu(a)}- q^{a_\nu(a)+1}
\over
t^{l_\mu(b)}-q^{a_\mu(b)+1}
}
=
{t^{\aaa_{k+1}-\aaa_{k+1}} -q^{\bbb_{k+1}-\bbb_k}
\over
t^{\aaa_{k+1}-\aaa_{\l+1}} -q^{\bbb_{\l}-\bbb_k}
}
=
{1
\over
t^{\aaa_{k+1}}}
{
\ol u_k-\ol x_{k+1}
\over
\ol u_k-\ol u_\l
}.
$$

\vskip -4.3in
\hfill $
\vcenter{\hbox{\includegraphics[width=1.7 in]{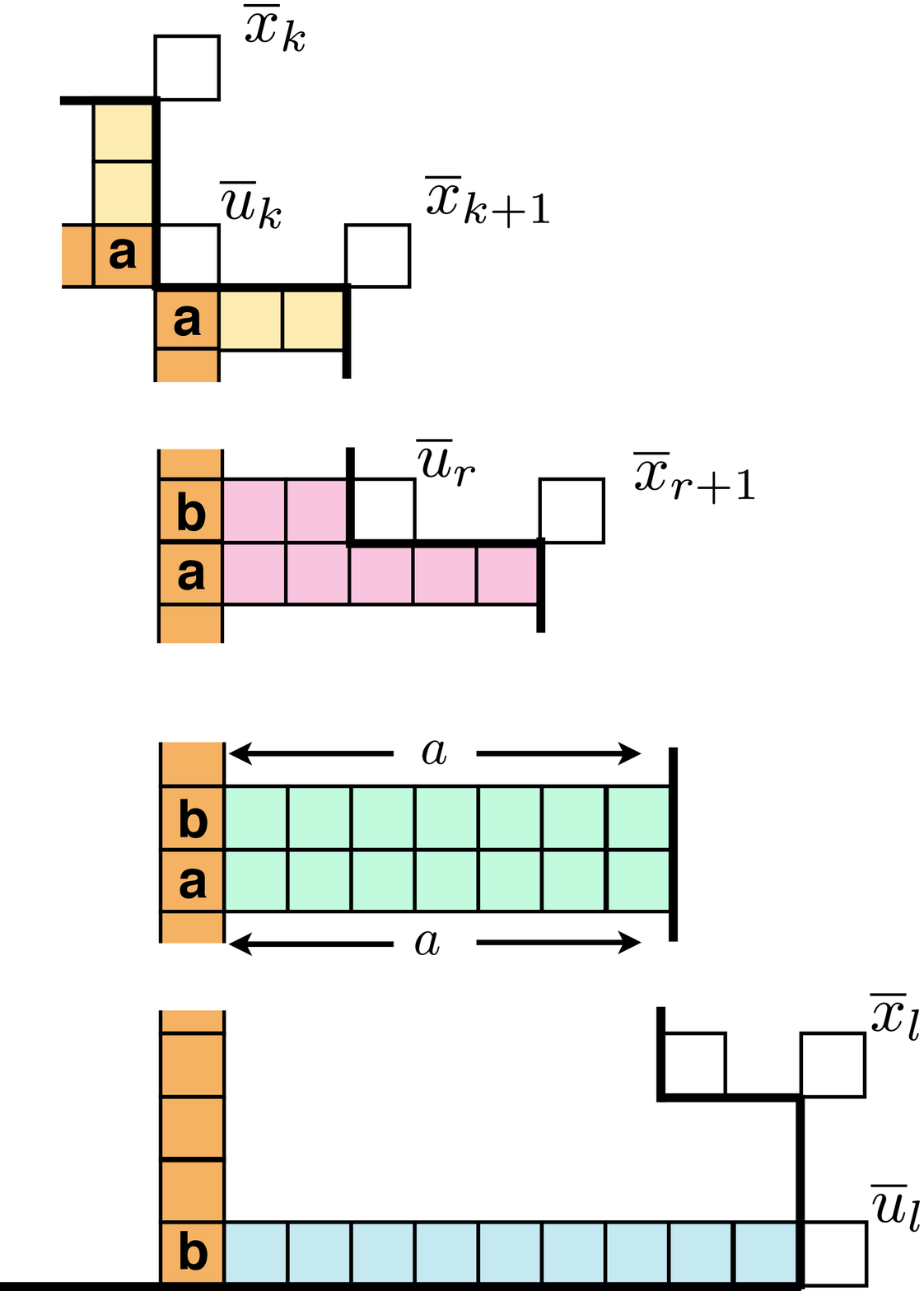}}}.
 $ 
\vskip 1.96in

Thus collecting all these ratios gives   
$$
 \prod_{s\in C_{\mu\nu}} 
{ 
t^{l_\nu(s)}-q^{a_\nu(s)+1}
\over 
t^{l_\mu(s)}-q^{a_\mu(s)+1}
}
\ses
{1
\over
t^{\aaa_{k+1}}}
{\prod_{i=k+1}^{\l}(\ol u_k-\ol x_i)
\over
{\prod_{i=k+1}^{l}(\ol u_k-\ol u_i)}
}
\eqno  3.14
$$
and we finally can see that the combination of 3.13 and  3.14 gives precisely the identity displayed in 3.9.
This completes our  proof of  3.6.
\sa

The most important consequence of the identity in 3.6 is 
the summation formula   \hbox{in 2.14 b)} which we are now in a position to derive with a minimum of efforts. The idea is to consider the rational function
$$
R(y)\ses
 {
\prod_{i=1}^\l (1-y \ol x_i)
\over
\prod_{i=0}^\l (1-y\ol u_i)
}
$$
and note that its partial fraction expansion may be written in the form
$$
R(y)\ses \sum_{j=0}^\l {A_j\over 1-y \ol u_j}
$$ 
with
$$
A_j=(1-y \ol u_j)R(y)\Big|_{y=1/\ol u_j} \bbu =  {
\prod_{i=1}^\l (1-  \ol x_i/\ol u_j)
\over
\prod_{i=0,i\neq j}^\l (1- \ol  u_i/ \ol u_j)
}
\ses \ol u_j\,  d_{\mu^ {(j)},\nu}.
 $$
Since  from 3.7 we derive that $\ol u_j={T_{\mu^ {(j)}}\over T_\nu}$ we see that 
the left hand side of 2.14 b) is none other than
$$
\sum_{j=0}^\l A_j \ol u_j^{k-1}\ses R(y)\Big|_{y^{k-1}}.
$$
But now the identity in 1.9 together with 3.8 gives
$$
 \eqalign{
 R(y)\Big|_{y^{k-1}}
 &\ses\OM\big[-y(  \ol x_1+\cdots +\ol x_\l\sms \ol u_0-\cdots -\ol u_\l)  \Big]
 \big|_{y^{k-1}}
 \cr
&\ses\OM\big[-y( MB_\mu-1 )  \big]\Big|_{y^{k-1}}
=h_{k-1}[-D_\mu]\ses (-1)^{k-1}e_{k-1}[D_\mu].
 \cr }
$$
This proves the first case of 2.14 b). The second case is immediate.
\sas

Here and in the following it will be convenient to set
$$
\Pi (a_1,a_2,\ldots ,a_n)\ses\nabla^{a_n} e_1 \nabla^{-a_n}\cdots \nabla^{a_2} e_1 \nabla^{-a_2}
\nabla^{a_1} e_1\nabla^{-a_1}.
\eqno 3.15
$$
As we have seen Theorem 1.1 assures that every operator $Q_{m,n}$ may be expressed as a non commutative polynomial in
the operators $Q_{m,1}$, thus the operators in 3.15 span 
the algebra $\CA$ generated by the operators $D_k$.
\sas

It turns out that the action of these operators have a remarkably beautiful Macdonald polynomial expansion,
a particular case of which may be stated as follows 
\sa

\noindent{\bol Theorem  3.1}

{\ita For any weak composition $a_1,a_2,\ldots ,a_n$  we have }
$$
 \Pi(a_1,a_2,\ldots ,a_n)\, {\bf 1}=\sum_{\mu\part n}\TH_\mu[X;q,t]
\sum_{T\in ST(\mu)}\prod_{i=2}^n{x_i^{1-a_i}\over 1-x_i}
\prod_{1\le i<j\le n}\OM[-Mx_j/x_i]\prod_{i=2}^n\big(1-x_i w_T(i)\big)\Big|_{S_T}
\eqno   3.16
$$
{\bol Proof}

 We may write
$$
\eqalign{
 \Pi(a_1,a_2,\ldots ,a_n)\, {\bf 1}
&=
 \sum_{\mu\part n}{\TH_\mu[X;q,t]\over w_\mu}
\ess \LL \Pi(a_1,a_2,\ldots ,a_n)\, {\bf 1}\scs \TH_\mu\RR_*
\cr
&=
\sum_{\mu\part n}{\TH_\mu[X;q,t]\over w_\mu}
\ess  \Pi(a_1,a_2,\ldots ,a_n)^{\perp^*}\TH_\mu.
\cr}
$$
Since the identity $P_{0,1}=-e_1$ gives
 $\ess
\nabla^a e_1 \nabla^{-a} = - \nabla^a P_{0,1} \nabla^{-a} ,
 \ess$
from 1.23 we derive that
$$
\big(\nabla^a e_1 \nabla^{-a}\big)^{\perp^*}\ses -\nabla^{-a }P_{0,-1} \nabla^{a}
\ses M\nabla^{-a }e_1^\perp \nabla^{a}.
$$
Thus
$$
 \Pi(a_1,a_2,\ldots ,a_n)^{\perp^*}\TH_\mu\ses 
M^n \nabla^{-a_1 }e_1^\perp \nabla^{a_1}\nabla^{-a_2 }e_1^\perp \nabla^{a_2}
\cdots \nabla^{-a_n }e_1^\perp \nabla^{a_n}\TH_\mu,
$$
and we can write 
$$
 \Pi(a_1,a_2,\ldots ,a_n)\, {\bf 1}
 =
\sum_{\mu\part n}{\TH_\mu[X;q,t]\over w_\mu}
\ess 
M^n \nabla^{-a_1 }e_1^\perp \nabla^{a_1}\nabla^{-a_2 }e_1^\perp \nabla^{a_2}
\cdots \nabla^{-a_n }e_1^\perp \nabla^{a_n}\TH_\mu.
\eqno 3.17
$$
To apply the operator $\nabla^{-a }e_1^\perp \nabla^{a}$ to $\TH_\mu$ we use 2.12 b) and
obtain, for any integer $a\ge 0$,
$$
\nabla^{-a }e_1^\perp \nabla^{a}\TH_\mu\ses \sum_{\nu\RA \mu}
c_{\mu,\nu}\big(\tttt{T_\mu\over T_\nu}\big)^a\TH_\nu.
$$
Starting from this identity, a straightforward induction argument  (carried out first in [9])  gives 
$$
 M^n\nabla^{-a_1 }e_1^\perp \nabla^{a_1}\nabla^{-a_2 }e_1^\perp \nabla^{a_2}
\cdots \nabla^{-a_n }e_1^\perp \nabla^{a_n}
\TH_\mu
\ses  M^n \sum_{T\in ST(\mu)}\prod_{k=2}^n w_T(k)^{a_k}c_{T^{(k)},T^{(k-1)}}
\eqno 3.18
$$
where $T^{(k)}$ denotes the tableau obtained from $T$ by removing all the entries larger than $k$ and
for notational convenience we have set
$$
c_{T^{(k)},T^{(k-1)}}\ses c_{shape(T^{(k)}),shape(T^{(k-1)})}.
$$
Next notice that a telescoping effect based on the 
 fact that $w_{T^{(1)}} =M$ yields the identity 
$$
{M\over w_\mu}\ses \prod_{k=2}^n
 {
 w_{T^{(k-1)}}
 \over 
 w_{T^{(k)}}
 }
 \bigsp \hbox{(with 
 $w_{T^{(k)}}=
 w_{
 shape(T^{(k)})}$),}
$$
which together with 2.13 written in the form
$$
Mc_{\mu,\nu}\ses {w_\mu\over w_\nu}d_{\mu,\nu}  
$$
allows us to carry out the following steps
$$
\eqalign{
 {M^n\over w_\mu}\nabla^{-a_1 }e_1^\perp \nabla^{a_1}\nabla^{-a_2 }e_1^\perp \nabla^{a_2}
\cdots \nabla^{-a_n }e_1^\perp \nabla^{a_n}
\TH_\mu
&\ses 
  {M \over w_\mu}   \sum_{T\in ST(\mu)}\prod_{k=2}^n w_T(k)^{a_k}Mc_{T^{(k)},T^{(k-1)}}
 \cr
&\ses 
 {M \over w_\mu} \sum_{T\in ST(\mu)}\prod_{k=2}^n w_T(k)^{a_k}
 {
 w_{T^{(k)}}
 \over 
 w_{T^{(k-1)}}
 }
d_{T^{(k)},T^{(k-1)}}
 \cr 
&\ses
\sum_{T\in ST(\mu)}\prod_{k=2}^n w_T(k)^{a_k}d_{T^{(k)},T^{(k-1)}}.
 \cr}
 \eqno 3.19
$$
The next move is to express $d_{T^{(k)},T^{(k-1)}}$ by means
of the formula given in   2.12, rewritten as follows
$$
\eqalign{
d_{\mu^{(k)},\nu}
&= {1\over \ol u_k} 
{\prod_{i=1}^\l\big(1-\tttt{\ol x_i\over \ol u_k}\big)\over  \prod_{i=0;i\neq k}^\l\big(1-\tttt{\ol  u_i\over \ol u_k}\big)}
\cr
&=
{z\over 1-z}
{\prod_{i=0}^\l\big(1-\tttt{\ol x_iz}\big)\over 
 \prod_{i=0}^\l\big(1-\tttt{\ol  u_iz}\big)}
 (1-\ol  u_kz)  \Big|_{z={1\over  \ol u_k}}
\cr
\hbox{(by 1.9)}
&=
{z\over 1-z}
\OM\big[
z(\ol u_0+ \cdots +\ol u_\l\sms \ol x_0-\cdots -\ol x_\l)
\big] (1-\ol  u_kz)\Big|_{z={1\over  \ol u_k}}
\cr
\hbox{(by 3.8)}
&=
{z\over 1-z}
\OM\big[
-zM B_\nu
\big] (1-\ol  u_kz)\Big|_{z={1\over  \ol u_k}}.
\cr}
\eqno 3.20
$$
 Since $\ul u_k=w_T^{(k)}$,
multiple substitutions  of 3.20  gives
$$
\eqalign{
{M^n\over w_\mu}
\nabla^{-a_1 }e_1^\perp \nabla^{a_1}&\nabla^{-a_2 }e_1^\perp \nabla^{a_2}
\cdots \nabla^{-a_n }e_1^\perp \nabla^{a_n}
\TH_\mu\ses
\cr
&=\sum_{T\in ST(\mu)}\prod_{k=2}^n w_T(k)^{a_k}
{z_k\over 1-z_k}\, \OM\big[-z_kMB_\nu \big](1-w_T(k)z_k)\Big|_{z_k={1\over  w_T(k)}}
\cr
&=\sum_{T\in ST(\mu)}\prod_{k=2}^n
{1\over z_k^{a_k-1}(1-z_k)}\, \prod_{1\le h<k}\OM\big[-Mz_k/z_h\big]\prod_{k=2}^n(1-w_T(k)z_k)\Big|_{z_k=1/w_T(k)}
\cr}
$$
which used in 3.18 gives 3.16 as desired.
 This completes our proof.

\sa

At this point it is  worthwhile also seeing what  tableaux expansions can  be obtained directly from  the $c_{\mu,\nu}$ formula given in 3.4. Recall that in  3.4  we have  depicted    the diagram
 of the partition  $\mu$ together with the $m(=5)$ cells 
that must be removed from  $\mu$ to obtain the partitions $\nu$ immediately preceding $\mu$ in the Young order. We  denoted there by $x_1,\ldots ,x_m$ their respective weights
 (from left to right). In addition to the removable cells, we have depicted the addable cells of $\mu$ that 
 are SW shifted by one unit, whose weights we denote by $u_0,u_1,\ldots ,u_m$. To the former sequence we have added the cell $(-1,-1)$  whose weight $1/qt$ we denoted by $x_0$. This given the $c_{\mu,\nu}$ formula  given in 3.4 can be more conveniently rewritten as
 
 \hsize=6.5in
 $$
 c_{\mu,\nu^{(i)}}={qt x_i\over M}
 {\prod_{j=0}^\l
 1-u_j/x_i
 \over 
 \prod_{j=1;j\neq i}^\l
 1-x_j/x_i
 }
\, =\,  \tttt{qt\over z M}\, (1-z/qt) {\prod_{j=0}^\l
 1-u_jz
 \over 
 \prod_{j=0}^\l
 1-x_jz
 }
(1-x_i z)\Big|_{z={1\over x_i}},
 \eqno 3.21
 $$
where $\nu^{(i)}$ is the partition obtained by removing from  $\mu$ the cell of weight $x_i$. Now to carry out the same sequence of steps that yielded the identity in 3.16 we start by using 3.5 and subject  the right hand side of 3.21 to the following 
successive transformations
$$
\eqalign{
c_{\mu,\nu^{(i)}}\, 
&=\,  \tttt{qt\over z M}\, (1-z/qt)
\OM\big[\tttt{z}(x_0+x_1+\cdots +x_l\sms u_0-u_1-\cdots -u_\l)
\big](1-x_iz)\Big|_{z={1\over x_i}}
\cr
&= 
 \tttt{qt\over z M}\,(1-z/qt)
 \OM
 \big[
\tttt{z\over qt}MB_\mu
 \big](1-x_iz)\Big|_{z={1\over x_i}}
 \cr
 (\hbox{Using $B_{\mu}=B_{\nu^{(i)}}+x_i$})
  &= 
 \tttt{qt\over z M}\,(1-z/qt)
 \OM
 \big[
\tttt{z\over qt}M B_{\nu^{(i)}}+\tttt{z\over qt}Mx_i 
 \big](1-x_iz)\Big|_{z={1\over x_i}}
 \cr
&= 
 \tttt{qt\over z M}\,(1-z/qt)
 \OM
 \big[
\tttt{z\over qt}M B_{\nu^{(i)}}
\big]\OM\big[zx_i-zx_i/t-zx_i/q+zx_i/qt  
 -x_iz\big] \Big|_{z={1\over x_i}}
 \cr
 &= 
 \tttt{qt\over z M}\,(1-z/qt)
 \OM
 \big[
\tttt{z\over qt}M B_{\nu^{(i)}}
\big]\OM\big[  -1/t-1/q+1/qt  
 \big] \Big|_{z={1\over x_i}}
 \cr
 &
\ses
 \tttt{qt\over z M}\,(1-z/qt)
 \tttt{(1-1/t)(1-1/q)\over 1-1/qt}
\OM
 \big[
\tttt{z\over qt}M B_{\nu^{(i)}} 
 \big]
 \ses
 \tttt{x_i }\, 
 \tttt{ 1-1/x_iqt )\over 1-1/qt}
\OM
 \big[
\tttt{M\over qt} B_{\nu^{(i)}}/x_i 
 \big].
 \cr}
$$
Using this identity    for $i\RA k$ and $x_i\RA w_T(k)$, 
3.18 becomes
$$
\eqalign{
 \nabla^{-a_1 }e_1^\perp \nabla^{a_1} & \nabla^{-a_2 }e_1^\perp \nabla^{a_2} 
\cdots \nabla^{-a_n }e_1^\perp \nabla^{a_n}
\TH_\mu
\cr
&
 \ses
  \sum_{T\in ST(\mu)}\prod_{k=2}^n w_T(k)^{a_k+1}
 \tttt{ 1-1/x_kqt )\over 1-1/qt}
\OM
 \Big[
\tttt{M\over qt} \sum_{1\le h\le k-1} w_T(h)/w_T(k) 
 \Big]
 \cr
&
 \ses
  \sum_{T\in ST(\mu)}\prod_{k=2}^n w_T(k)^{a_k }
 \tttt{ 1-w_T(k) qt   )\over 1-qt }
\OM
 \Big[
\tttt{M\over qt} \sum_{1\le h\le k-1} w_T(h)/w_T(k) 
 \Big]
 \cr 
 & \ses
  \sum_{T\in ST(\mu)}\prod_{k=2}^n w_T(k)^{a_k }
 \tttt{ 1-w_T(k) qt   )\over 1-qt }
\prod_{1\le h<k \le n}
\OM
 \Big[
\tttt{M\over qt}\tttt{  w_T(h)\over w_T(k) }
 \Big]
 \cr 
}
$$
which  in turn, substituted in 3.17, gives
$$
 {1\over M^n}\Pi(a_1,a_2,\ldots ,a_n)\, {\bf 1}
 \ses
\sum_{\mu\part n}{\TH_\mu[X;q,t]\over w_\mu}
\sum_{T\in ST(\mu)}\prod_{k=2}^n w_T(k)^{a_k }
 \tttt{ 1-w_T(k) qt   )\over 1-qt }
\prod_{1\le h<k \le n}
\OM
 \Big[
\tttt{M\over qt}\tttt{  w_T(h)\over w_T(k) }
 \Big].
$$
This proves Theorem I.3 stated in the introduction.
\sas

The above results give only a glimpse of the computational  power of   $d_{\mu,\nu}$ and $c_{\mu,\nu}$.
 Note that in the original paper [16]
Macdonald derived Pieri formulas  
for multiplication of his polynomials 
by $e_r[X]$ and $h_r[\tttt{1-t\over 1-q}X]$,
for any  $r\ge 1$. Nevertheless, for many years we have had  Pieri coefficients corresponding  only  to multiplication (or skewing) of $\TH_\mu[X;q,t]$  by  $e_1$. 
Actually there is a simple reason for this. In passing from   $P_\mu[X;q,t]$  to $\TH_\mu[X;q,t]$
the Macdonald formulas yield Pieri coefficients corresponding to multiplication (or skewing) by $h_r[\tttt{X\over 1-t}]$ and
$e_r[\tttt{X\over 1-q}]$ and it is only when $r=1$ that the
Macdonald  formulas yield Pieri coefficients corresponding to multiplication by an elementary or the homogeneous symmetric function.
 
Actually it turns  out
   that multiplication (or skewing) by $e_1$ is not at all as
limited as it may  appear  on the surface. Indeed, in  a recent (2012) paper [5],  Bergeron-Haiman, guided by Hilbert Scheme considerations, discovered that multiplication by
$e_r[\tttt{X\over M }]$ and skewing by $h_r[X]$ for any $r\ge 1$
may be recursively expressed in terms of the $d_{\mu,\nu}$ and $c_{\mu,\nu}$. Since this development is very closely related to our operators $Q_{m,-1}$ and is conducive to a wide variety of standard tableaux expansions, we will terminate this 
section and the paper with a brief presentation of some applications of the  Bergeron-Haiman identities.

For a given $k\ge 1$ let us set
$$
a)\ess\ess h_k^\perp \TH_\mu[X;q,t]= \sum_{\nu\con_k \mu}\ck \TH_\nu [X;q,t]
\ess\ess\ess\ess\ess\ess
b)\ess\ess  e_k[\tttt{X\over M }]  \TH_\nu[X;q,t]=\sum_{\mu\noc_k \nu}\dk \TH_\mu [X;q,t]
\eqno 3.22
$$ 
where ``$\nu\con_k \mu$'' means that $\nu$ is contained in $\mu$ (as Ferrers diagrams) and $\mu/\nu$ has $k$ lattice cells. The symbol ``$\mu\noc_k \nu$'' is analogously defined.
It follows from the orthogonality of the basis $\{\TH_\mu[X;q,t]\}_\mu$ that in this case 3.2 becomes
$$
c_{\mu,\nu}^{(k)}\ses \tttt{w_\mu\over w_\nu}d_{\mu,\nu}^{(k)}
\eqno 3.23
$$ 
This given, the Bergeron-Haiman identities  may be 
stated as follows
\sas

\noindent{\bol Theorem 3.2}

{\ita For any $k\ge 1$ and $\mu\part n$ we have}
$$
c_{\mu,\nu}^{(k+1)}\, =\,   {1\over B_{\mu/\nu}}
\sum_{\nu\con_1\aaa \con_{ k }\mu
}c_{\mu,\aaa}^{(k)}c_{\aaa,\nu}^{(1)}T_\aaa/T_\nu
\ess\ess\ess\ess\big(\hbox{with $B_{\mu/\nu}=B_ \mu-B_ \nu $}\big)
\eqno 3.24
$$

In spite of the cumbersome denominator, this is a remarkably 
rapid  way of obtaining the $c_{\mu,\nu}^{(k)}$ in terms 
of our ${c_{\mu,\nu}}'s$. In particular, 3.24   yields by far the fastest algorithm for computing the $\TH_\mu[X;q,t]$ to this date.
With any efficient Symbolic Manipulation software, in a few seconds we may obtain the monomial 
expansion of  $\TH_\mu[X;q,t]$ for $\mu\part n$  with $n$ as large as $18$. In fact, the coefficients $M_{\la,\mu}(q,t)$ in the expansion
$$ 
\TH_\mu[X;q,t]\ses \sum_{\la\part n}m_\la[X]M_{\la,\mu}(q,t)
\ess\ess\ess\ess( \hbox{for $\mu\part n $})
\eqno 3.25
$$
can be simply calculated by means of the identities
$$
M_{\la,\mu}(q,t)\ses \LL h_\la\scs \TH_\mu\RR
\ses h_\la^\perp \TH_\mu.
\eqno 3.26
$$

The most surprising development concerning this 2012 discovery is that it is an immediate consequence  of a simple identity obtained   in the 1990's!   The latter identity is best 
stated, as  expressed in [6],  in terms of the generating functions
$$
\CT(u) P[X] =P[X +u]=\sum_{k\ge 0}u^k h_k^\perp P[X]
\bigsp
D(z)= \sum_{r}z^r D_r = \OM[-zX]\CT(\tttt{M\over z}),
\eqno 3.27
$$
where for any expression $E$ we set
$$
\CT(E) P[X] =P[X +E].
$$
In fact,  the commutativity of the two translations $\CT(u)$
and $\CT(\tttt{M\over z})$ immediately yields
$$
\CT(u)D(z)= \OM[-z(X+u)]\CT(\tttt{M\over z})\CT(u)=
(1-uz)D(z)\CT(u)
$$
or better
$$
\big[
D(z),\CT(u)
\big]= uz D(z)\CT(u).
\eqno 3.28
$$
This brings us in position to  include here a  proof of 
Theorem 3.2.  Equating the coefficients of $u^{k+1}z^0$
in 3.28 gives (using 1.11 $(i)$ and $(iii)$)
$$
\big[
(I-M\Delta_{e_1}),h_{k+1}^\perp
\big]= 
\big[
D_0,h_{k+1}^\perp
\big]= 
  D_{-1}h_k^\perp=M\nabla^{-1}e_1^\perp \nabla h_k^\perp
$$
which may be rewritten as
$$
\big[h_{k+1}^\perp,\Delta_{e_1}\big]=\nabla^{-1}e_1^\perp \nabla h_k^\perp
\eqno 3.29
$$
and 3.24 immediately follows from the definition in 3.22 a)
by applying both sides of 3.29 to $\TH_\mu$.
\sas

It will be instructive to see in a more explicit way how the recursion in 3.24 expresses   $c_{\mu,\nu}^{(k)}$ as a polynomial in the ${c_{\mu,\nu}}'s$. This is best carried out as a standard tableaux expansion. To begin, given a standard tableau $T$ with $n$ cells let $T^{(k)}$ as before denote the 
standard tableau obtained from $T$ by removing all the entries 
greater than $k$, in particular we may set $T=T^{(n)}$.
Also set
$$
c_{T^{(k)},T^{(k-1)}}=c_{shape(T^{(k)}),shape(T^{(k-1)})}
\ess\ess\ess\hbox{and}\ess\ess\ess
B_{T^{(k)}/T^{(h)}} =B_{shape(T^{(k)})}-B_{shape(T^{(h)})}
$$
and define for $1\le m\le n$
$$
\PI_m(T;q,t)=\prod_{k=0}^{m-1}
\left(
{c_{T^{(n-k)},T^{(n-k-1)}}
\over 
B_{T^{(n)}/T^{(n-k-1)}}
}wt(T^{(n-k)})
\right)
\eqno 3.30
$$ 
where ``$wt(T^{(s)})$'' denotes the weight of the cell we must remove from $T^{(s)}$  to obtain  $T^{(s-1)}$.

Using this notation we can state the following consequence
of the recursion in 3.24
\supereject

\noindent{\bol Proposition 3.1}
$$
h_m^\perp\TH_\mu[X;q,t]
\, =\,
\sum_{T\in\ST(\mu)}
\Pi_m(T;q,t)
\TH_{shape(T^{(n-m)} )}[X;q,t]
\eqno 3.31
$$
From this we can immediately derive the following
corollary.

\sas

\noindent{\bol Theorem 3.3}

{\ita If $\mu\part n$ and $\la=(\la_1,\la_2,\ldots ,\la_k)\part n$ 
then setting $p_i=\la_1+\la_2+\cdots +\la_{i-1}$, the coefficients $M_{\la,\mu}(q,t)$ in 3.25 yielding the monomial expansion  
of $\TH_\mu[X;q,t]$ are given by  the following standard tableaux sums}
$$
M_{\la,\mu}(q,t)\, =\,
\sum_{T\in ST(\mu)}\prod_{i=1}^k
\Pi_{\la_i}(T^{(n-p_i)};q,t).
\eqno 3.32
$$
This leads to the following efficient way to obtain the expansion of any symmetric function in terms of the modified Macdonald 
basis.

For any given  $F[X]\in \La^{=n}$, start with the expansion 
$$
F[X]\ses \sum_{\mu\part n}{\TH_\mu[X,q,t]\over w_\mu}
\LL \TH_\mu\scs F\RR_*,
$$
then use 3.25 and obtain
$$
F[X]\ses \sum_{\mu\part n}{\TH_\mu[X,q,t]\over w_\mu}
\sum_{\la\part n}
\LL m_\la\scs F\RR_* M_{\la,\mu}(q,t).
$$

With precomputed $ M_{\la,\mu}(q,t)$ according to 3.32
this permits experimenting with various combinatorial questions in Macdonald polynomial theory working with symmetric functions of substantially larger  degrees than was possible before the  Bergeron-Haiman 2012 discovery.
\sap

\centerline{\bol Bibliography}

\item{[1]} F. Bergeron and A. M. Garsia  {\ita Science Fiction and Macdonald's Polynomials}.
 CRM Proceedings  $\&$ Lecture Notes, American Mathematical
  Society , 22:, (1999),  1--52.
\sas

\item{[2]} F. Bergeron, A. M. Garsia, M. Haiman and G. Tesler, {\ita Identities and positivity conjectures for some remarkable operators in the theory of symmetric functions}, Methods Appl. Anal. 6, (1999), 363--420.
\sas

\item{[3]}
 F. Bergeron, A. M. Garsia ,  E. Leven and G. Xin,   
{\ita Some remarkable new Plethystic Operators in the Theory of Macdonald Polynomials},
arXiv:1405.0316 [math.CO]
\sas

\item{[4]}
 F. Bergeron, A. M. Garsia ,  E. Leven and G. Xin,   
{\ita A Compositional  (km,kn)-Shuffle Conjecture},
arXiv:1404.4616 [math.CO]
\sas

\item{[5]}
F. Bergeron, M. Haiman,{\ita  Tableaux Formulas for Macdonald Polynomials}, 
Special edition in honour of Christophe Reutenauer 60 birthday, International Journal of Algebra and Computation, Volume 23, Issue 4, (2013), pp 833--852.
\sas

\item{[6]} 
A.~M. Garsia, M.~Haiman, and G.~Tesler.
 {\ita Explicit plethystic formulas for Macdonald $(q,t)$-Kostka
  coefficients}.
 S\'eminaire Lotharingien de Combinatoire [electronic only],
  42:\penalty0 B42m, 1999.
\sas

\item{[7]}
A. M. Garsia and  M. Haiman,
{\ita A graded representation model for Macdonald's polynomials}, Proc.  Nat. Acad. Sci. U.S.A {\bol 90} no. 8,
 (1993), 3607--3610.
\sas

\item{[8]}
A. M. Garsia and  M. Haiman,
{\ita Some Natural Bigraded $S_n$-Modules and q,t-Kostka Coefficients. (1996) (electronic)
The Foata Festschrift}, http://www.combinatorics.org/Volume 3/volume 3 2.html\#R24
\sas

\item{[9]} 
A. M. Garsia, J. Haglund  and G.~Xin, 
{\ita Constant term methods in the theory of Tesler matrices and Macdonald polynomial operators}, Annals of Combinatorics, DOI: 10.1007/s00026-013-0213-6, (2013).
\sas

\item{[10]} 
A. M. Garsia and G. Tesler, {\ita Plethystic formulas for Macdonald q,t-Kostka coefficients}, Adv. Math. {\bol  123}
no. 2, (1996), 144--222.
\sas

\item{[11]} 
A. M. Garsia, M. Zabrocki   and G.~Xin, 
 {\ita HallÐLittlewood Operators in the Theory of Parking Functions and Diagonal Harmonics}, 
 Int Math Res Notices (2011)
doi: 10.1093/imrn/rnr060
\sas

\item{[12]} 
E. Gorsky and A.~Negut.
 {\ita Refined knot invariants and Hilbert schemes}.
 arXiv preprint arXiv:1304.3328, 2013.
\sas

\item{[13]} 
J. Haglund, M.~Haiman, N.~Loehr, J.~B. Remmel, and A.~Ulyanov.
 {\ita A combinatorial formula for the character of the diagonal
  coinvariants}.
 Duke J. Math., 126:\penalty0, (2005), 195--232.
\sas

\item{[14]} 
J.~Haglund, J.~Morse, and M.~Zabrocki.
 {\ita A compositional refinement of the shuffle conjecture specifying touch
  points of the dyck path.}
 Canadian J. Math, 64:\penalty0, (2012), 822--844.
\sas

\item{[15]} 
T.~Hikita.
 {\ita Affine springer fibers of type a and combinatorics of diagonal
  coinvariants.}
 arXiv preprint arXiv:1203.5878, 2012.
\sas

\item{[16]}
I. G. Macdonald, {\ita A new class of symmetric functions}, Actes du 20e S\' eminaire
Lotharingien, Publ. I.R.M.A. Strasbourg 372/SÐ20, (1988), 131--171.
\sas

\item{[17]}
I. G. Macdonald , {\ita Symmetric functions and Hall polynomials}, Oxford Mathematical Monographs, Oxford Science Publications, Oxford University Press, New York, 1995.
\sas

\item{[18]} 
A.~Negut,
{\ita  The shuffle algebra revisited},
 arXiv preprint arXiv:1209.3349, 2012.
\sas

\item{[19]} 
O.~Schiffmann and E.~Vasserot.
 {\ita The elliptical Hall algebra, Cherednik Hecke algebras and Macdonald
  polynomials}.
 Compos. Math., 147.1:\penalty0, (2011), 188--234.
\sas

\item{[20]} 
O.~Schiffmann and E.~Vasserot.
 {\ita The elliptical Hall algebra and the equivariant K-theory of the
  Hilbert scheme of ${A}^2$}.
 Duke J. Math., 162.2:, (2013), 279--366.
\sas

\item{[21]} 
O. Schiffmann, {\ita On the Hall algebra of an elliptic curve, II}, Preprint (2005),
arXiv:math/0508553.
\sas

\item{[22]} 
R. Stanley, {\ita Some Combinatorial Properties of Jack Polynomials}, Advances in Mathematics 
{\bol  77},  (1989), 76--115.

\end